\RequirePackage{fix-cm}
\documentclass[smallextended]{svjour3}       
\smartqed  
\usepackage{graphicx}

\usepackage{graphicx,epsfig}
\usepackage{amsmath,amssymb}
\usepackage{exscale}

\smartqed
\usepackage{soul}
\usepackage{amssymb}
\usepackage{mathptmx}
\usepackage{color}
\usepackage{float}
\usepackage{subfigure}

\usepackage{algorithm,algpseudocode}
\usepackage{tikz}
\usepackage{pgfplots}

\usepackage{amsfonts}
\usepackage{amssymb}
\usepackage{amsmath}

\usepackage{algorithm}
\usepackage{algorithmicx}
\usepackage{algpseudocode}

\usepackage{multirow}
\usepackage{longtable}

\usepackage{rotating, tabularx}

\usepackage{enumitem}

\newcommand{\sgn}{\operatorname{sgn}}

\newcommand{\hypo}{\operatorname{hypo}}
\newcommand{\graph}{\operatorname{graph}}
\newcommand{\epi}{\operatorname{epi}}

\newtheorem{Illustrative Example}{Illustrative Example}

\begin{document}


\title{	 On radially epidifferentiable functions and regularity conditions in nonconvex optimization\thanks{This study is supported by The Scientific and Technological Research Council of Turkiye (TUBITAK) under the Grant No. 217M487.}}


\titlerunning{On radially epidifferentiable functions}        

\author{Gulcin Dinc Yalcin        \and
        Refail Kasimbeyli 
}


\institute{G. Dinc Yalcin and R. Kasimbeyli  \at
	Faculty of Engineering,
	Department of Industrial Engineering,
	Eskisehir Technical University, Iki Eylul Campus, Eskisehir 26555, Turkey \\
	Tel.: +90-222-3213550/6445\\
	\email{gdinc@eskisehir.edu.tr}\\
	Corresponding Author: Refail Kasimbeyli \at
	Affiliation: Department of Industrial Engineering,
	Eskisehir Technical University, Eskisehir 26555, Turkey \\ \email{rkasimbeyli@eskisehir.edu.tr}          
}

\date{Received: date / Accepted: date}

\maketitle

\begin{abstract}
	In this paper we study the radial epiderivative notion for nonconvex functions, which extends the (classical) directional derivative concept. The paper presents new definition and new  properties for this notion and establishes relationships between the radial epiderivative, the Clarke's directional derivative, the Rockafellar's subderivative and the directional derivative. The radial epiderivative notion is used to establish new regularity conditions without convexity conditions. The paper analyzes necessary and sufficient conditions for global optimums in nonconvex optimization via the generalized derivatives studied in this paper. We establish a necessary and sufficient condition for a descent direction for radially epidifferentiable nonconvex functions. The paper  presents explicit formulations for computing the weak subgradients in terms of the radial epiderivatives and vice versa, which are very important from point of view of developing solution methods for finding global optimums in nonsmooth and nonconvex optimization. All the properties and theorems presented in this paper, are illustrated and interpreted on examples.
\keywords {Generalized derivative \and Rockafellar's subderivative \and Clarke's directional derivative \and radial epiderivative \and weak subgradient \and nonconvex optimization}
\subclass{90C26 \and 90C30 \and 90C46 \and 90C56.}
\end{abstract}

\section{Introduction} \label{intro}

Historically, the derivative was and remains to be a powerful tool in analyzing optimization problems and plays an important role in developing methods for characterizing and computing optimal solutions. 
If the function $f:\mathbb{R}^n \rightarrow \overline{\mathbb{R}}$ is convex and differentiable at $ \overline{x}$ with the gardient vector   $\nabla f(\overline{x}),$  then the graph of the affine function  $a(x)=f(\overline{x}) + \langle \nabla f(\overline{x}), x-\overline{x} \rangle$ becomes a supporting hyperplane to $epi f$ at $(\overline{x},f(\overline{x}))$ and for every $x \in dom f$ one has  
\begin{equation}\label{gradineq}
f(x) -  f(\overline{x}) \geq \langle \nabla f(\overline{x}), x-\overline{x} \rangle,
\end{equation}
or $epi f \subset epi a$ and $f(\overline{x}) = a(\overline{x}).$

For a variety of reasons, the differentiability condition may be considered as too restrictive in applictions. This situation has led to the development of the theory of generalized
differentiation and hence the development of the theory of nonsmooth analysis (see e.g. \cite{AbbasiKT2021,Aubin1981,BorweinI1996,Clarke1975,EkelandL1976,Floresbazan2003,KabganiL2022,KabganiS2022,Kasimbeyli2009,LaraK2021,LinhP2013,Luc1991,MichelP1992,Mordukhovich1984,Penot1978,Rockafellar1980}).

At the first step of such a generalization, the limit $t \rightarrow 0$ was relaxed by changing it to the one-sided form $t \downarrow 0,$ which will be called in this paper, a directional derivative $f^{\prime}(\overline{x};h)$ of  $f$ at $\overline{x}$ in a direction $h \in \mathbb{R}^n$  and defined as
\begin{equation}\label{dirderdefintro}
f^{\prime}(\overline{x};h) =  \lim_{ t \downarrow 0} \frac{f(\overline{x}+ th)- f(\overline{x})}{t},
\end{equation}
if the limit exists.

The relation \eqref{gradineq} was extended to a nondifferentiable case by T. Rockafellar in \cite{Rockafellar1963}, where the subdifferential $\partial f(\overline{x})$ was defined as a set of subgradient vectors $v \in \mathbb{R}^n:$
\begin{equation}\label{rocksubdif}
\partial f(\overline{x}) = \{v \in \mathbb{R}^n : f(x) - f(\overline{x}) \geq \langle v, x-\overline{x} \rangle, \forall x \in \mathbb{R}^n \}.
\end{equation}
This notion successfully extends the affine support relations given above for the gradient vector $\nabla f(\overline{x})$ and the affine function $a(x).$ The subdifferentiability of $f$ at $\overline{x}$ means the existence of a vector  $v \in \mathbb{R}^n$ such that the  hyperplane 
$$
H(v,-1) = \{(x,y) : \langle (v,-1), (x-\overline{x}, y - f(\overline{x}))\rangle =0 \}
$$  
with normal vector $(v,-1),$ is a supporting hyperplane to the epigraph of $f$ at $(\overline{x},f(\overline{x})):$ 
\begin{equation}\label{episubsetaffine}
epi f \subset H^-(v,-1) = \{(x,y) : \langle (v,-1), (x-\overline{x}, y - f(\overline{x}))\rangle  \leq 0 \}.
\end{equation}
The global nature of this relation leads to the following necessary and sufficient condition for unconstrained global optimality of $\overline{x}$: 
\begin{equation}\label{optcondzerosubgrad}
0 \in \partial f(\overline{x}) \Leftrightarrow f(x) \geq f(\overline{x}), \forall x\in \mathbb{R}^n.
\end{equation}
For the problem of minimizing a convex function $f(x)$ subject to $x \in S, $ where $S \subset \mathbb{R}^n$ is a convex, closed set, the following necessary and sufficient condition for  global optimality was formulated by Rockafellar \cite{Rockafellar1970}  and Pshenichnyi \cite{Pshenichnyi1971}: $\overline{x} \in S$ is a global optimal solution to this problem if and only if there exists a subgradient $v \in \partial f(\overline{x})$ such that 
\begin{equation}\label{pshrockoptcond}
\langle v, x - \overline{x} \rangle \geq 0, \forall x\in S.
\end{equation}

Subgradients in the convex case found widespread applications in optimization (see e.g. \cite{EkelandT1974,Rockafellar1970,RockafellarW2009}). In the convex case, the subdifferential set $\partial f(\overline{x})$ can be characterized with respect to the directional derivatives of $f$ at $\overline{x}:$ 
\begin{equation}\label{subdifrockdirder1}
\partial f(\overline{x}) =  \{v \in \mathbb{R}^n : f^{\prime}(\overline{x};h) \geq \langle v,h \rangle \mbox{  for all  } h \in \mathbb{R}^n \},
\end{equation}
\begin{equation}\label{subdifrockdirder2}
f^{\prime}(\overline{x};h) = \max \{\langle v,h \rangle : v \in \partial f(\overline{x}) \} \mbox{  for all  } h \in \mathbb{R}^n. 
\end{equation}

A major contribution to the nonsmooth and nonconvex analysis was made by F.H. Clarke.  In \cite{Clarke1973}, Clarke introduced a generalized directional derivative concept $f^{\circ}$ and  showed how the definition of subdifferential $\partial f$ can be extended to arbitrary lower semicontinuous, locally Lipschitz functions defined on Banach spaces. Clarke introduced the notion of the generalized subdifferential $\partial^{\circ}f(\overline{x})$ as a set of subgradient vectors $v \in \mathbb{R}^n$ with
\begin{equation*}\label{clarkesubdifintro}
\partial ^{\circ}f(\overline{x}) = \{v \in \mathbb{R}^n : f^{\circ}(\overline{x};h) \geq \langle v,h \rangle, \forall h \in \mathbb{R}^n \}.
\end{equation*}
The main property of the Clarke stationary points is given in the following theorem.

\begin{theorem} {\normalfont\cite[Proposition 2.3.2, p. 38]{Clarke1983}}\label{clarke-stationary-point}
	If $f$ attains a local minimum or maximum at $\overline{x}$ then $0 \in \partial^{\circ}f(\overline{x}).$	
\end{theorem}

Clarke introduced the regularity notion which plays an important role in nonsmooth analysis.

\begin{definition}{\normalfont\cite[Definition 2.3.4, p. 39]{Clarke1983}}\label{regularfintro}
	$f$ is said to be regular at $x$ provided that the classical directional derivative $f^{\prime}(x;h)$ exists and $f^{\circ}(x;h) = f^{\prime}(x;h) $ for all $h.$
\end{definition}

The following regularity theorem was given in \cite{Clarke1983}.

\begin{theorem} {\normalfont\cite[Proposition 2.3.6 (b), p. 40]{Clarke1983}}\label{regularftheorem1}
	Let $f$ be Lipschitz near $x.$ If $f$ is convex, then $f$ is regular at $x.$	
\end{theorem}

The regularity property for the Clarke's directional derivative was proved under the convexity condition, and hence the relation with the directional derivative, is not applicable in a nonconvex case. 
On the other hand, although Clarke extended the subgradient notion to Lipschitz continuous functions, there was nothing analogous for the case of general $f.$ Rockafellar made a serious contribution to fill this gap, by introducing a subderivative function $df(\overline{x};\cdot)$ \cite{Rockafellar1979,Rockafellar1980}. Bu using the subderivative notion, Rockafellar established a new regularity condition for the Clarke's directional derivative in the Lipschitzian case (see Theorem \ref{subderregular} in Section \ref{regularity}). The extensions of the derivative notion given by Clarke and Rockafellar, have made a huge contribution to the nonsmooth analysis. However, even under the regularity conditions, the stationary points $\overline{x}$ with $f^{\circ}(\overline{x};h) \geq f^{\circ}(\overline{x};0) = 0$  or $\text{d}f(\overline{x};h) \geq \text{d}f(\overline{x};0) =0  $ for all $h,$ characterize only local minimum points of $f.$

The analysis given above demonstrates that, for an optimization problem with nonconvex and nondifferentiable functions, to formulate conditions guaranteeing the existence of supporting surfaces similar to \eqref{gradineq} and/or \eqref{episubsetaffine}, conditions for global optimality similar to \eqref{optcondzerosubgrad} and/or \eqref{pshrockoptcond}, or characterization relations similar to \eqref{subdifrockdirder1} and/or \eqref{subdifrockdirder2}, are not easy tasks and require additional assumptions and new ideas. The characterization of a global minimum by using the derivatives and/or generalized derivatives, still remains to be one of the main problems in the mathematical programming. It can be expected that such a characterization may also help to develop solution methods for escaping from local minima. One of the main purposes of this paper, is to study tools, which allow to analyze these problems.

Two of such tools studied in this paper, are the  radial epiderivative and the weak subdifferential concepts, both introduced earlier by R.Kasimbeyli.

By replacing affine functions in the definition \eqref{rocksubdif} of $\partial f$ by  superlinear functions, R. Kasimbeyli weakened the affine support relations by changing them to a conical support ones to the epigraph of $f$ and introduced the concept of the weak subdifferential $\partial^wf$  in his dissertation \cite{Gasimov1992} (see also \cite{AzimovG1999,AzimovG2002,DincyalcinK2020mmor,DincyalcinK2020opt}), as a set of weak subgradients $(v,c) \in \mathbb{R}^n \times \mathbb{R}_+ $ with
\begin{equation*}\label{WSeqintro}
f(x)\geq f(\overline{x})+\langle v,x-\overline{x} \rangle -c \|x-\overline{x}\|, \quad \forall x \in \mathbb{R}^n.
\end{equation*}
The existence of a weak subgradient at $\overline{x},$ corresponds to the existence of a supporting conical surface to $\epi f$ at $(\overline{x},f(\overline{x}))$ and allows by this way to handle applications not fitting within the domain of convexity. By using the conical supporting philosophy developed in  \cite{Gasimov1992}, Kasimbeyli extended the Hahn-Banach separation theorem to a nonconvex case \cite{Kasimbeyli2010,KasimbeyliK2019} which has played an important role in analyzing nonconvex optimization problems and developing solution methods for them \cite{DincyalcinK2020mmor,DincyalcinK2020opt,Gasimov2002,Kasimbeyli2013,KasimbeyliM2009,KasimbeyliM2011}. With the help of the weak subgradients, the global optimality condition \eqref{optcondzerosubgrad} was extended to a nonconvex case (see \cite[Remark 2.2]{DincyalcinK2020opt}):
\begin{equation}\label{optcondzeroweaksubgrad}
(0,0) \in \partial^w f(\overline{x}) \Leftrightarrow f(x) \geq f(\overline{x}), \forall x\in \mathbb{R}^n.
\end{equation}

Recently, Dinc Yalcin and Kasimbeyli developed a new  weak subgradients based global solution method for nonconvex box-constrained optimization problems in \cite{DincyalcinK2020opt}, where the approximate computing method for the weak subgradients via the directional derivatives, was also suggested.

The radial epiderivative concept was first proposed by F. Flores-Bazan in \cite{Floresbazan2001} (see also  \cite{Floresbazan2003}). In this paper, we use the definition of this concept given by Kasimbeyli in \cite{Kasimbeyli2009} (in a slightly different setting, for both set-valued maps and real-valued functions). By using the radial epiderivative $f^r(\overline{x};h)$ of a function $f: \mathbb{R}^n \rightarrow \overline{\mathbb{R}}$ at a point $ \overline{x}$ in a direction $h\in \mathbb{R}^n,$ a necessary and sufficient conditions for the global minimum $\overline{x}$ of a real-valued nonconvex function $f$ is given in \cite{Floresbazan2001} and \cite{Kasimbeyli2009}:
\begin{equation}\label{optcondzeroradepider}
f^r(\overline{x};h) \geq f^r(\overline{x};0) = 0, \forall h\in \mathbb{R}^n \Leftrightarrow f(x) \geq f(\overline{x}), \forall x\in \mathbb{R}^n.
\end{equation}
In this paper, this condition is used to establish a descent direction for a nonconvex funcion.


An extension of the optimaity condition \eqref{pshrockoptcond} to a nonconvex case, was established by Kasimbeyli and Mammadov in terms of the weak subgradients \cite[Theorem 4]{KasimbeyliM2011}. Kasimbeyli and Mammadov proved under some mild conditions that, the radial epiderivative and the directional derivative  can be represented as a support function of the  weak subdifferential set, and by this way the characterization relations similar to \eqref{subdifrockdirder1} and/or \eqref{subdifrockdirder2}, in nonconvex case were established by using the directional derivatives, radial epiderivatives and the weak subdifferentials \cite[Theorems 4.5 and 4.6]{KasimbeyliM2009}. 

Although all the above mentioned theorems and properties were established by using the definition of the radial epiderivative given in \eqref{radepider}, it would be interesting to have a definition of this concept, formulated in terms of the conventional limit relation of the Newton quotient $ \frac{f(\overline{x}+ tu)- f(\overline{x})}{t}. $ In this paper, we give such a formulation for the radial epiderivative and prove that $f^r(\overline{x};\cdot)$ is a lower semicontinuous and lower Lipschitz function. We study the relations between the radial epiderivatives and the  Clarke's directional derivative, the subderivatives and the directional derivatives and establish new regularity conditions. All these relations allow to connect the contributions made by the radial epiderivative concept and the huge world of the nonsmooth analysis due to the classical directional derivative, the Clarke directional derivative and the Rockafellar's subderivative. By taking into account the above mentioned relations earlier established between the radial epiderivatives and the weak subgradients, and the weak subgradient based globally solution method in nonconvex analysis, the methods for computing the weak subgradients becomes an important issue. This paper presents new constructive theorems which enable to compute the weak subgradients via the radial epiderivatives and vice versa. Computing the weak subgradients via the radial epiderivatives has an advantage over the computation by using the directional derivatives. The reason is that the directional derivative allows to estimate the local behaviour of the given function, while the radial epiderivative characterizes the global behaviour. The paper presents a comprehensive analysis on the proved theorems and established properties by using illustrative examples.


\looseness-1 The rest of the paper is organized as follows. The main definitions are given in Section \ref{pre}. Section \ref{secradialepider} presents a new formulation and some properties of the radial epiderivative. The relations between the radial epiderivatives and the directional derivative, Clarke's directional derivative and the Rockafellar's subderivative are presented in Section \ref{regularity}. In Section \ref{computingweaksubgr}, we establish conditions guaranteeing the existence of a weak subgradient for radially epidifferentiable functions and present formulas for approximate computing the radial epiderivatives in terms of weak subgradients and vice versa. These conditions are proved by considering both the Euclidean ($\ell_2$)  and the $\ell_1-$norms separately. Section \ref{optcond} presents a theorem on a globally descent direction and discusses a condition on the global optimality for nonconvex functions. Finally, Section \ref{conclusion} draws some conclusions from the paper.


\section{Preliminaries} \label{pre}

We begin this section by first recalling the definitions and some important properties of the main concepts used in this paper.

\begin{definition}\label{allcones}
	Let $S$ be a nonempty subset of a real normed space $(\mathbb{X},\|\cdot \|)$ and $\overline{x} \in S$ be a given element. The closed radial cone $R(S; \overline{x})$ of $S$ at $\overline{x}$ is the set of all $w \in \mathbb{X} $ such that there are sequences $\lambda_n > 0$ and  $(x_n)_{n \in \mathbb{N}} \subset S$ with $\lim_{n \rightarrow +\infty} \lambda_n (x_n- \overline{x}) = w.$ In other words,
		$$
		R(S; \overline{x})=cl(cone(S-\overline{x})),
		$$
		where \textit{cone} denotes the conic hull of a set, which is the smallest cone containing $S-\overline{x}.$
\end{definition}

Now, we recall the definitions of the generalized derivatives which will be investigated in this paper. We begin with the Clarke's directional derivative.

\begin{definition} {\normalfont \cite{Clarke1983}} \label{clarkedirdersubdifintro}
	Let $\mathbb{X}$ be a Banach space, let $f: \mathbb{X} \rightarrow \mathbb{R} $ be a locally Lipschitz function and let $\overline{x} \in \mathbb{X}$ and $h \in \mathbb{X}$ be given elements. The Clarke directional derivative $f^{\circ}(\overline{x};h)$ of $f$ at $\overline{x}$ in the direction $h,$ is defined by
	\begin{equation*}\label{clarkedirderintro}
	f^{\circ}(\overline{x};h)=\limsup_{t \downarrow 0, y \rightarrow \overline{x}} \frac{f(y+th) - f(y)}{t}. 
	\end{equation*}
\end{definition}

It was proved that, the Clarke directional derivative $f^{\circ}(\overline{x};\cdot)$ is upper semicontinuous, positively homogeneous, sublinear and Lipschitz function \cite[Proposition 2.1.1, p. 25]{Clarke1983}. The convexity of $f^{\circ}(\overline{x};\cdot)$ was proved by Rockafellar \cite[Theorem 1]{Rockafellar1979b}.

Now, we recall the definiton of Rockafellar's subderivative \cite[Definition 8.1, p.299]{RockafellarW2009} (see also \cite{Rockafellar1979,Rockafellar1980}).

\begin{definition}\label{subderdef}
	For a function $f: \mathbb{R}^n  \rightarrow \mathbb{R} $ and  a point $\overline{x} \in \mathbb{R}^n $ with $f(\overline{x})$ finite, the subderivative $df(\overline{x};h)$ of function $f$ at $\overline{x}$ in a direction $h \in \mathbb{R}^n, $ is defined by
	\begin{equation*}\label{subderint}
	\text{d}f(\overline{x};h)= \liminf_{t\downarrow 0, u \rightarrow h} \frac{f(\overline{x}+tu) - f(\overline{x})}{t}. 
	\end{equation*}
\end{definition}

\begin{remark}
	The definition of the subderivative given in Definition \ref{subderdef}, is more specifically the lower subderivative of $f$ given in \cite{Rockafellar1980}, where the corresponding upper subderivative is defined with ``lim sup'' in place of ``lim inf''. Moreover, T. Rockafellar used the term ``the subderivative of $f$ at $\overline{x}$ for $h$'' for the subderivative $df(\overline{x};h)$ \cite[Definition 8.1, p.299]{RockafellarW2009}.
\end{remark}

\begin{definition}\label{radepiderdef}
The radial epiderivative $f^r(\overline{x};h)$ of a function $f: \mathbb{R}^n \rightarrow \overline{\mathbb{R}}$ at a point $ \overline{x}$ in a direction $h\in \mathbb{R}^n,$ is defined through the radial cone $R(epi f; (\overline{x},f(\overline{x})))$ to the epigraph $epi f$ of $f$ at $(\overline{x},f(\overline{x}))$ such that 
\begin{equation}\label{radepider}
epi f^r(\overline{x};\cdot) = R(epi f; (\overline{x},f(\overline{x}))).
\end{equation}
 In the case when the radial epiderivative $f^{r}(\overline{x};h)$ exists and finite for every $h,$ we will say that $f$ is radially epidifferentiable at $\overline{x}.$
\end{definition}
The radial epiderivative is probably the first derivative concept which extends the global affine suport relations \eqref{gradineq} and \eqref{episubsetaffine} to a nonconvex case by using a global conical supporting surface to the epigraph of a function under consideration.

\begin{remark}\label{radepiderremark}
	In the original definition of the radial epiderivative given in \cite{Kasimbeyli2009}, the notation $D_r F( \overline{x}, \overline{y})(\cdot)$ was used for this notion, which was defined for a set-valued map $F,$ where $\overline{y} \in F(\overline{x})$. Since in this paper we consider real-valued functions, we use the notation $f^{r}(\overline{x};\cdot).$ Such a notation is similar to those used for the directional derivative, the Clarke's and the Rockafellar's derivatives, and by this way we aim to use unified notations for all the generalized derivatives considered in this paper. 
\end{remark}

Now we recall the existence condition for the radial epiderivative proved in \cite{Kasimbeyli2009}.

\begin{theorem}{\normalfont \cite[Theorem 3.2]{Kasimbeyli2009}\label{rk3.2}}
	Let $(\mathbb{X}, \|.\|_\mathbb{X})$ be a real normed space,  $\mathbb{S}$ be a non-empty subset of $\mathbb{X}, \quad \overline{x} \in \mathbb{S}$  and let $f:\mathbb{S} \rightarrow \mathbb{R} \cup \{+\infty\}$  be a proper function. Assume that there exist functions $g_1,g_2: \mathbb{X} \rightarrow \mathbb{R}$ with $epi(g_1) \subset R(epi(f); (\overline{x}, f(\overline{x}))) \subset epi(g_2).$ Then the radial epiderivative $f^r(\overline {x};\cdot)$ is given as
	\begin{equation} \label{radial20}
	f^r(\overline {x};h)= \inf\{ y\in \mathbb{R}: (h,y) \in R(epi(f);(\overline{x}, f(\overline{x}))) \},  \forall h \in \mathbb{X}.
	\end{equation}
\end{theorem}

\begin{lemma} {\normalfont \cite[Lemma 3.7]{KasimbeyliM2009}\label{radposhom}}
	Let $f:\mathbb{X}\rightarrow \mathbb{R}$ be a single-valued function having radial epiderivative
	$f^r(\overline {x};\cdot)$ given by \eqref{radial20}.  Then the radial
	epiderivative $f^r(\overline {x};\cdot)$ is a positively homogeneous function.
\end{lemma}


The generalized derivatives given by Clarke and by Rockafellar, were used to define the corresponding generalized subgradients. These generalized subgradients were defined as the normal vectors of supporting  hyperplanes to the epigraphs of the corresponding derivatives. In difference to these concepts, the classical subgradient of a function in the convex analysis introduced by Rockafellar, was defined as the normal vector of the supporting hyperplane to the epigraph of the function under consideration. It is remarkable that in nonconvex analysis, this property was kept probably only in the definition of the weak subgradient, which strongly fits this property.

\begin{remark}\label{conicsupport}
	It follows from definition of the weak subgradient given in Section \ref{intro} that, the pair $(v,c) \in \mathbb{R}^n \times \mathbb{R}_+$ is a weak subgradient of $f$ at $\overline{x} \in \mathbb{R}^n,$ if there exists a continuous (superlinear) concave function 
	\begin{equation*}
	g(x)= f(\overline{x})+\langle v,x-\overline{x} \rangle -c \|x-\overline{x}\|,
	\end{equation*}
	such that $g(x) \leq f(x),$ for all $x \in \mathbb{R}^n$ and $g(\overline{x}) = f(\overline{x}).$ Then clearly the set $\hypo(g)=\{(x,\alpha) \in \mathbb{R}^n \times \mathbb{R} : g(x) \geq \alpha \}$ is a closed convex cone in $\mathbb{R}^n \times \mathbb{R}$ with vertex at $(\overline{x}, f(\overline{x})),$ and
	$$
	\epi(f) \subset \epi(g), \quad cl(\epi(f)) \cap \graph(g) \neq \emptyset .
	$$
\end{remark}


The above analysis shows that the class of weakly subdifferentiable functions is essentially larger than the class of subdifferentiable functions, see e.g.  \cite[Theorems 3.1, 3.2, Corollary 3.1]{AzimovG1999}, \cite[Theorem 1]{AzimovG2002}
\cite[Lemma 2.8]{KasimbeyliM2009}, \cite[Theorem 3]{DincyalcinK2020mmor}, \cite[Theorem 2.3]{DincyalcinK2020opt}.
The following theorem explines some classes of the weakly subdifferentiable functions, which will be used in the next sections.

\begin{theorem}{\normalfont \cite[Lemma 2.7]{KasimbeyliM2009}} \label{wsubdiffunc}
	Let $S$ be a nonempty subset of a real normed space $(\mathbb{X},\|\cdot \|)$. Let $f: \mathbb{S} \rightarrow (- \infty , + \infty]$ be a given function.  
	If $f$ is a positively homogeneous function bounded from below on some neighborhood of $0_{\mathbb{R}^n},$ then $f$ is weakly subdifferentiable at $0_{\mathbb{R}^n}.$
\end{theorem}



\section{Properties of radial epiderivatives}\label{secradialepider}

This section presents new properties of the radial epiderivative and some illustrative examples. We begin with the lower semicontinuity property for the radially epidifferentiable functions.


\begin{theorem} \label{radepiderlsc}
	Let $(\mathbb{X}, \|.\|_\mathbb{X})$ be real normed space, $\mathbb{S}$ be a non-empty subset of $ \mathbb{X},$ and let $f:\mathbb{S} \rightarrow \mathbb{R} \cup \{+\infty\}$  be a proper function. If $f$ is radially epidifferentiable at $\overline{x} \in \mathbb{X}$ then $f$ is lower semicontinuous at $\overline{x}.$ 
\end{theorem}	
\begin{proof}
	let $f:S \rightarrow \mathbb{R}$ be a proper function, radially epidifferentiable at $\overline{x}.$  We need to show that $\liminf_{ x \rightarrow \overline{x}} f(x) \geq f(\overline{x}).$ Assume to the contrary that this is not true: $\liminf_{ x \rightarrow \overline{x}} f(x) < f(\overline{x}).$ Then,  radial cone	$R(epi(f);(\overline{x}, f(\overline{x})))$ must contain a vertical line passing through the points $(\overline{x}, f(\overline{x})) $ and  $(\overline{x}, \liminf_{ x \rightarrow \overline{x}} f(x))$ with  
	\begin{equation*}
	\inf\{ y\in \mathbb{R}: (0,y) \in R(epi(f);(\overline{x}, f(\overline{x}))) \} = -\infty
	\end{equation*}
	which contradicts the hypothesis  that 	$f$ is radially epidifferentiable at $\overline{x},$ and hence the proof is completed.
	$ \Box $ 
\end{proof}


\begin{remark}
	Note that the inverse of Theorem \ref{radepiderlsc} is not true. For example $f(x) = -\sqrt{|x|}$ is (lower semi) continuous at $x=0$ but not radially epidifferentiable there.
\end{remark}


The following proposition gives an equivalent representation for the radial epiderivative via a limit concept. Note that a similar expression was also given by F. Flores-Bazan in terms of the lower epiderivative (see \cite[Corollary 3.4]{Floresbazan2003}).


\begin{proposition}\label{rews}
	Let $(\mathbb{X}, \|.\|_\mathbb{X})$  be a real normed space and let $\overline{x} \in \mathbb{X}$ be a given element. Assume that function $f:\mathbb{X} \rightarrow \mathbb{R} $ is radially epidifferentiable at $\overline{x}.$ Then the radial epiderivative $f^{r}(\overline{x},\cdot)$ can equivalently be defined as follows:
	
	\begin{equation} \label{radepilimdef}
	f^{r}(\overline{x};h) = \inf_{t > 0} \liminf_{ u \rightarrow h} \frac{f(\overline{x}+ tu)- f(\overline{x})}{t}
	\end{equation}
	for all $h \in \mathbb{X}.$ 
\end{proposition}


\begin{proof} Let $\overline{x} \in \mathbb{X}$ and $\overline{y} = f(\overline{x}.)$
	By the definition of the radial epiderivative, we have:
	
	\begin{eqnarray*}
		&&R(\epi(f);(\overline{x}, \overline{y}))= \epi(f^{r}(\overline{x};\cdot)) \\
		&=& \{ (x,y) \in \mathbb{X} \times \mathbb{R} : \exists  \lambda_n >  0,   (x_n ,y_n) \in \epi(f),  \lim_{n \rightarrow \infty} \lambda_n ((x_n, y_n)-(\overline{x} , \overline{y}))=(x,y) \} \\ 
		&=&\{(x,y) \in \mathbb{X} \times \mathbb{R} : \exists  \lambda_n >  0,   (x_n ,y_n) \in \epi(f),\\  
		&&\lim_{n \rightarrow \infty} \lambda_n (x_n-\overline{x})=x, \lim_{n \rightarrow \infty} \lambda_n (y_n- \overline{y})=y  \}.
	\end{eqnarray*}
	
	By Theorem \ref{rk3.2},
	\begin{eqnarray*}
		f^{r}(\overline{x};h) = \inf \{y: \lambda_n > 0, (x_n ,y_n) \in \epi(f),   &h&=\lim_{n \rightarrow \infty} \lambda_n (x_n-\overline{x}), \qquad \qquad \\
		&y&=\lim_{n \rightarrow \infty} \lambda_n (y_n- \overline{y})  \}.	
	\end{eqnarray*}
	The last equality can be written in the following form:
	\begin{equation*}
	f^{r}(\overline{x};h) = \inf \{ y :  \lambda_n >  0,   h =\lim_{n \rightarrow \infty} \lambda_n (x_n-\overline{x}),
	y=\lim_{n \rightarrow \infty} \lambda_n (f(x_n)- \overline{y})  \}.	
	\end{equation*}
	By setting $x_n= x_n -\overline{x} + \overline{x},$ we deduce:
	\begin{eqnarray*}
		f^{r}(\overline{x};h)	= \inf \Big\{ &y& :  \lambda_n >  0,  \\
		&h&=\lim_{n \rightarrow \infty} \lambda_n (x_n-\overline{x}), 
		y=\lim_{n \rightarrow \infty} \lambda_n (f(x_n -\overline{x} + \overline{x})- f(\overline{x}))  \Big\}	
	\end{eqnarray*}
	or
	\begin{eqnarray*}
		f^{r}(\overline{x};h) = \inf \Big\{ &y& :  \lambda_n >  0,     \\
		&h&=\lim_{n \rightarrow \infty} \lambda_n (x_n-\overline{x}), 
		y=\lim_{n \rightarrow \infty} \frac{f(\frac{\lambda_n (x_n -\overline{x} )}{\lambda_n} + \overline{x} )- f(\overline{x})}{1 / \lambda_n}  \Big\}.
	\end{eqnarray*}
	Letting $u_n= \lambda_n (x_n -\overline{x}),$ for every  $n=1,2, \ldots,$ we obtain:
	\begin{equation*}
	f^{r}(\overline{x};h)= \inf \Big\{y : \lambda_n >  0, \mbox{  for every  } n=1,2, \ldots,
	y=\liminf_{n \rightarrow \infty,  u_n \rightarrow h} \frac{f(\overline{x} +\frac{u_n}{\lambda_n} )- f(\overline{x})}{1 / \lambda_n}  \Big\}.
	\end{equation*}
	Now by setting $t_n = 1 / \lambda_n,$ for every $n=1,2, \ldots,$ we can rewrite the last relation as follows:
	\begin{eqnarray*}
		f^{r}(\overline{x};h) = \inf \Big\{ y :  t_n >  0,  \mbox{  for every  } n=1,2, \ldots,   y=\liminf_{u \rightarrow h} \frac{f(\overline{x} +t_n u)- f(\overline{x})}{t_n}\Big\},
	\end{eqnarray*}
	which completes the proof. $\Box$
\end{proof}


\begin{theorem} \label{radepiderlowerlip}
	Let $(\mathbb{X}, \|.\|_\mathbb{X})$ be a real normed space and  $f:\mathbb{X} \rightarrow \mathbb{R}$ be a
	proper function finite at $x=\overline{x}$.  If $f$ is lower Lipschitz at $\overline{x},$ that is there exists a positive constant $L$ such that 
	\begin{equation}\label{fislowerlip}
	f(x)-f(\overline{x}) \geq  -L \|x-\overline{x}\|
	\quad \mbox{ for all } x \in \mathbb{X},
	\end{equation}
	then $f$ is radially epidifferentiable at $\overline{x}.$ If $\mathbb{X} =\mathbb{R}^n$ then this condition is also necessary.
\end{theorem}	

\begin{proof}
	Assume that $f$ is lower Lipschitz at $\overline{x}:$ there exists a positive constant $L$ such that \eqref{fislowerlip} is satisfied for all $x \in \mathbb{X}.$ Take an arbitrary element $h \in \mathbb{X}$ and evaluate the expression
	\begin{equation} \label{radepilimdef1}
	\inf_{t > 0} \lim_{ u \rightarrow h} \frac{f(\overline{x}+ tu)- f(\overline{x})}{t}.
	\end{equation}
	By \eqref{fislowerlip} we obtain:
	\begin{equation*}
	\inf_{t > 0} \lim_{ u \rightarrow h} \frac{f(\overline{x}+ tu)- f(\overline{x})}{t} \geq \inf_{t > 0} \lim_{ u \rightarrow h} \frac{-tL\|u\|}{t} =-L\|h\|,
	\end{equation*} 
	which means that the expression \eqref{radepilimdef1} has a finite value, and hence we deduce by \eqref{radepilimdef} that, $f$ has a finite radial epiderivative at $\overline{x}$ in every direction $h.$
	
	Now assume that $X =\mathbb{R}^n$ and that $f$ has a finite radial epiderivative at $\overline{x}$ in every direction $h \in \mathbb{R}^n$. Show that 
	there exists a positive constant $L$ such that \eqref{fislowerlip} is satisfied for every $x \in \mathbb{R}^n.$ Assume to the contrary that this is not true. Let $\{x_n\} \subset \mathbb{R}^n$ and $\{L_n\}$ be sequences  with $L_n \rightarrow +\infty$  and $\|x_n-\overline{x}\| >0$ such that
	\begin{equation*}
	f(x_n)-f(\overline{x}) <  -L_n \|x_n-\overline{x}\|
	\quad \mbox{ for all } n =1,2, \ldots.
	\end{equation*}
	Let $u_n =  \frac{x_n-\overline{x}}{\|x_n-\overline{x}\|}.$ Without loss of generality, assume that  $u_n \rightarrow h$ with $\|h\|=1$ and let $t_n = \|x_n-\overline{x}\| >0$ for all $n.$
	Then, we obtain:
	\begin{equation*} 
	\lim_{ n \rightarrow \infty} \frac{f(\overline{x}+ t_nu_n)- f(\overline{x})}{t_n} <  \lim_{ n \rightarrow \infty} \frac{-t_n L_n\|u_n\|}{t_n} =-\infty,
	\end{equation*}
	which contradicts to the assumption that $f$ is radially epidifferentiable at $\overline{x}.$ 
	$\Box$
\end{proof}


\begin{remark}
	The lower Lipschitz concept, was called calmness, or calm from below in \cite[Chapter 8, Section F, p.322]{RockafellarW2009}.  
\end{remark}


Now consider some examples and demonstrate the properties of the radial epiderivative.

\begin{figure}[h]
	\scalebox{0.39}{\includegraphics{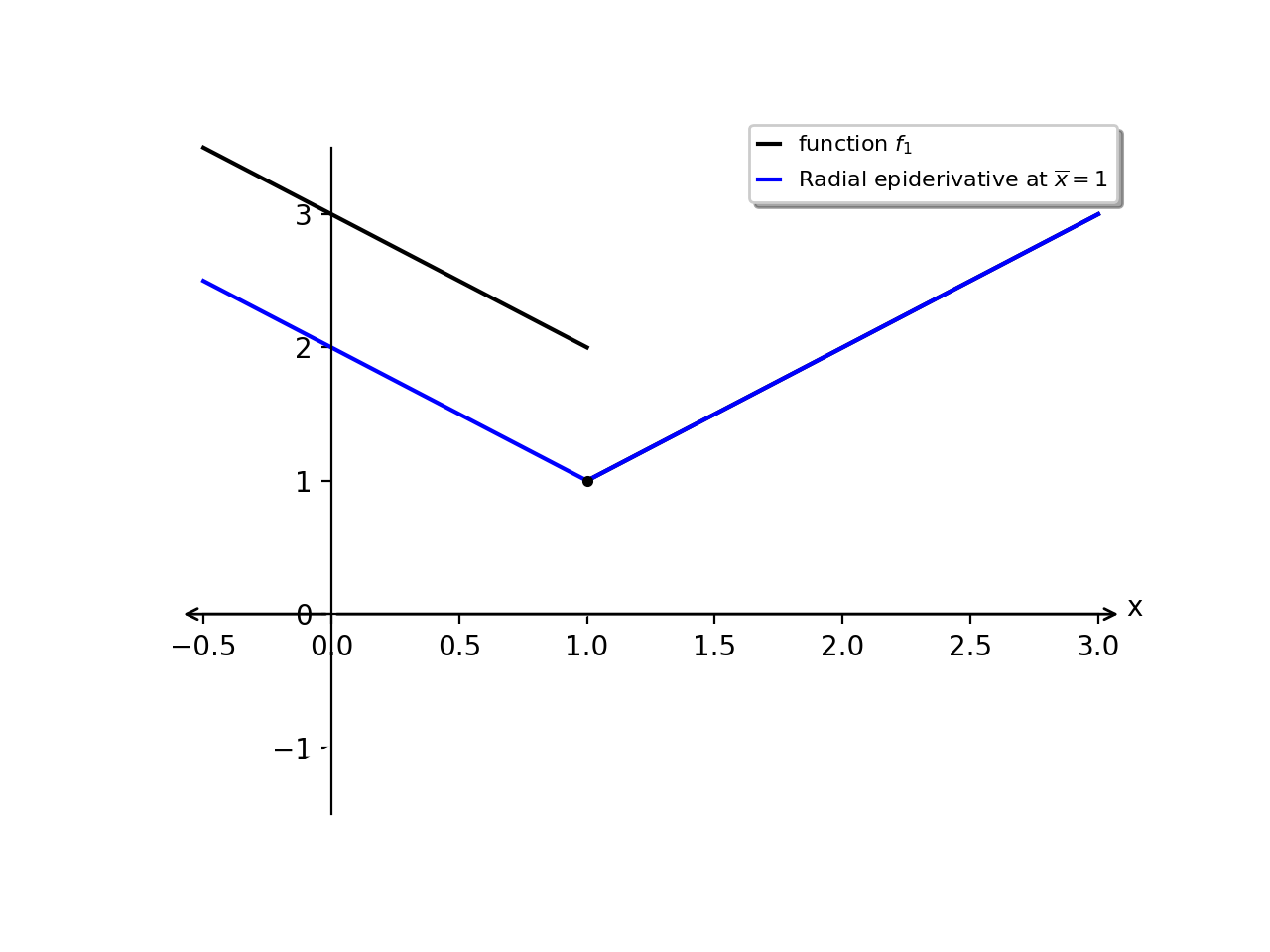}}
	\scalebox{0.39}{\includegraphics{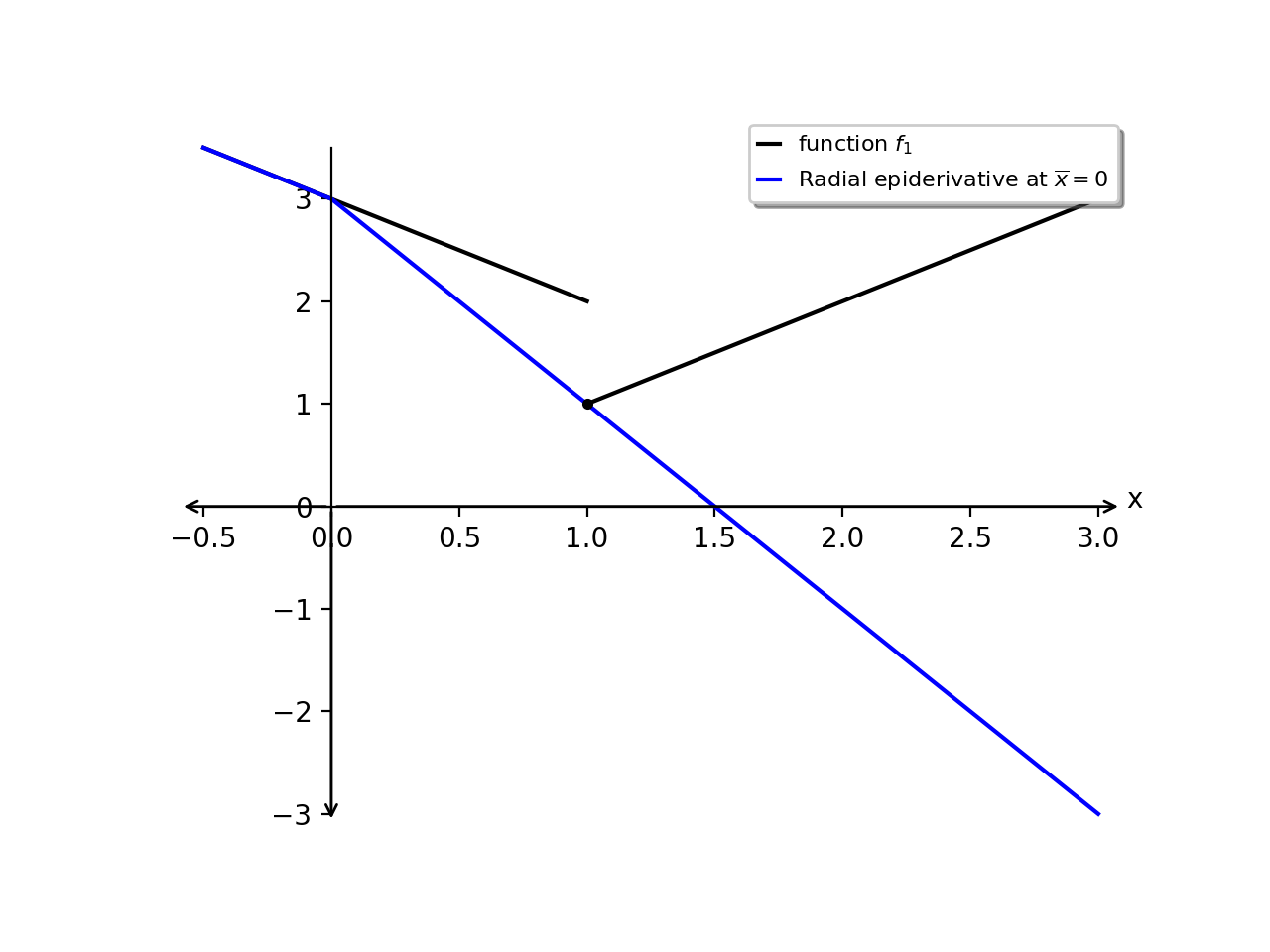}}
	\caption{The graph of radial epiderivative of function $f_1$ at $\overline{x}=0$ (left) and $\overline{x}=1$ (right)}
	\label{figure:f_1_RE}
\end{figure}

\begin{example} \label{example1}
	Let 
	\begin{equation*}\label{funcex1}
	f_1(x) =
	\left\{
	\begin{array}{ll}
	-x+3    & \mbox{ if } x < 1, \\
	x & \mbox{ if } x \geq 1.
	\end{array}
	\right.
	\end{equation*}
	
	The function $f_1$ is defined and lower semicontinuous everywhere on $\mathbb{R}.$ 
	This function is not continuous at $\overline{x} =1$ and does not satisfy the
	Lipschitz condition there. It is just lower Lipschitz at $\overline{x} =1.$ \\ 
	
	It follows from Theorems \ref{rk3.2}  and  \ref{radepiderlowerlip}  that $f_1$ has a radial epiderivative $f^r_1(\overline{x};\cdot)$ at every point $\overline{x} \in \mathbb{R}.$ By using Proposition \ref{rews}, we obtain (see also Figure \ref{figure:f_1_RE}):

	\begin{equation*}\label{radepider2ex1}
	f^r_1(\overline{x};h) =
	\left\{
	\begin{array}{llll}
	-h  & \mbox{ if } \overline{x} \leq 1, h < 0, \\
	\frac{(\overline{x} - 2)h}{1 - \overline{x}}    & \mbox{ if } \overline{x} < 1, h > 0, \\
	h & \mbox{ if } \overline{x} = 1, h>0, \\
	h & \mbox{ if } \overline{x} > 1, h\in \mathbb{R}. \\
	\end{array}
	\right.
	\end{equation*}

\end{example}

\begin{example} \label{example2}
	
	Let 
	\begin{equation*}\label{funcex2}
	f_2(x) =
	\left\{
	\begin{array}{ll}
	-x+3    & \mbox{ if } x \leq 1, \\
	x & \mbox{ if } x > 1.
	\end{array}
	\right.
	\end{equation*}
	Clearly, it follows from Theorem \ref{rk3.2} that $f_2$ has a radial epiderivative $f^r_2(\overline{x};h) = f^r_1(\overline{x};h)$ at every point $\overline{x} \neq 1, $ but  at $\overline{x} =1,$ where the function is not lower semicontinuous (and not lower Lipschitz), we have (see also Theorems \ref{radepiderlsc} and \ref{radepiderlowerlip}):
	$$
	f^r_2(1;h) =
	\left\{
	\begin{array}{ll}
	-h    & \mbox{ if } h \leq 0, \\
	-\infty & \mbox{ if } h > 0.
	\end{array}
	\right.
	$$
	
\end{example}


\begin{example} \label{example3}
	Let 
	\begin{equation*}\label{funcex3}
	f_3(x) =
	\left\{
	\begin{array}{ll}
	4|x+1|    & \mbox{ if } x \leq 0, \\
	|x-1| + 3 & \mbox{ if } x > 0.
	\end{array}
	\right.
	\end{equation*}

	It follows from Theorem \ref{rk3.2} that $f_3$ has a radial epiderivative $f^r_3(\overline{x};\cdot)$ at every point $\overline{x} \in \mathbb{R}.$ Again by applying Proposition \ref{rews}, we obtain (see Figure \ref{figure:f_3_RE} for illustrations):


	\begin{figure}
		\scalebox{0.39}{\includegraphics{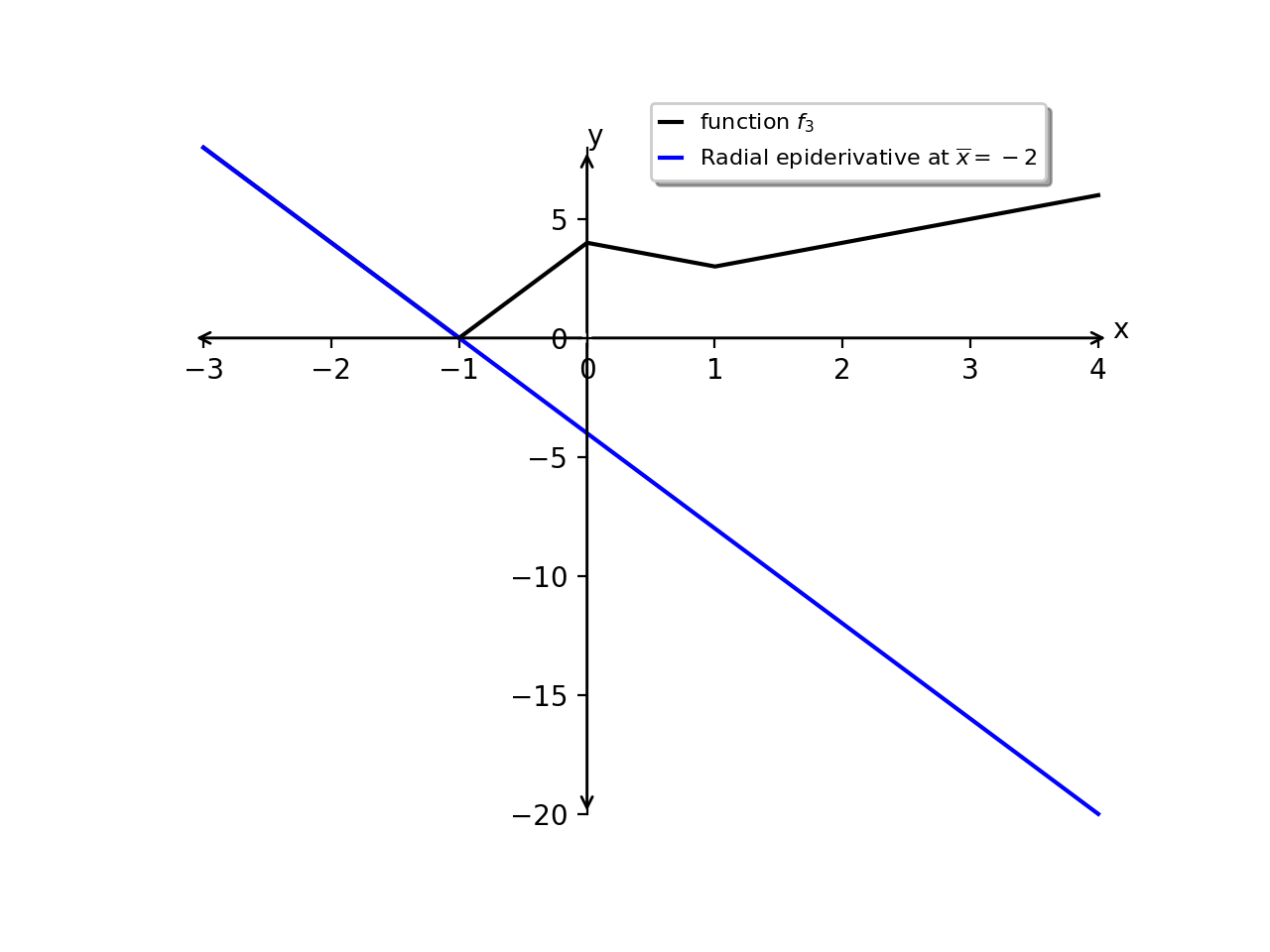}}
		\scalebox{0.39}{\includegraphics{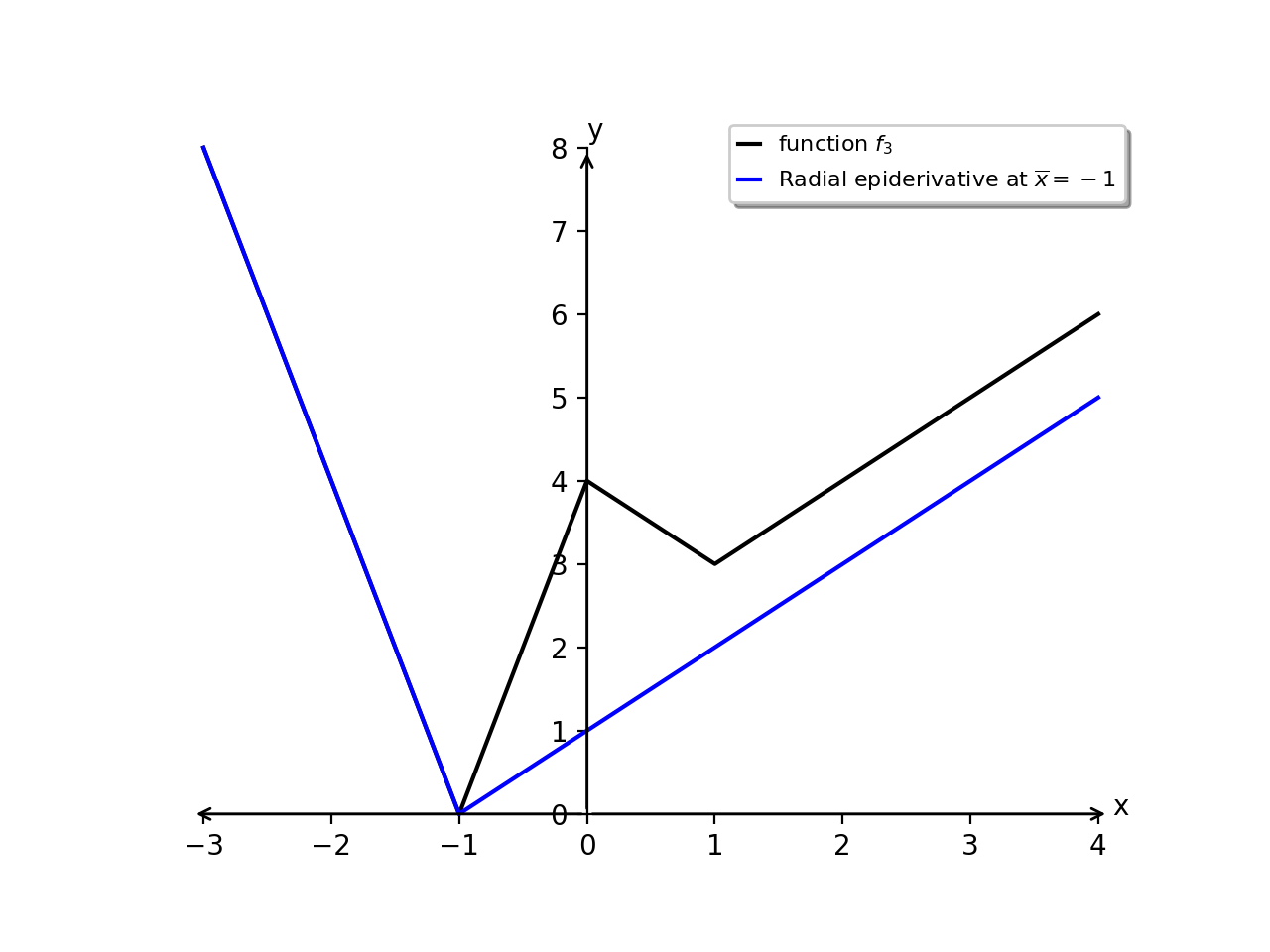}}
		\scalebox{0.39}{\includegraphics{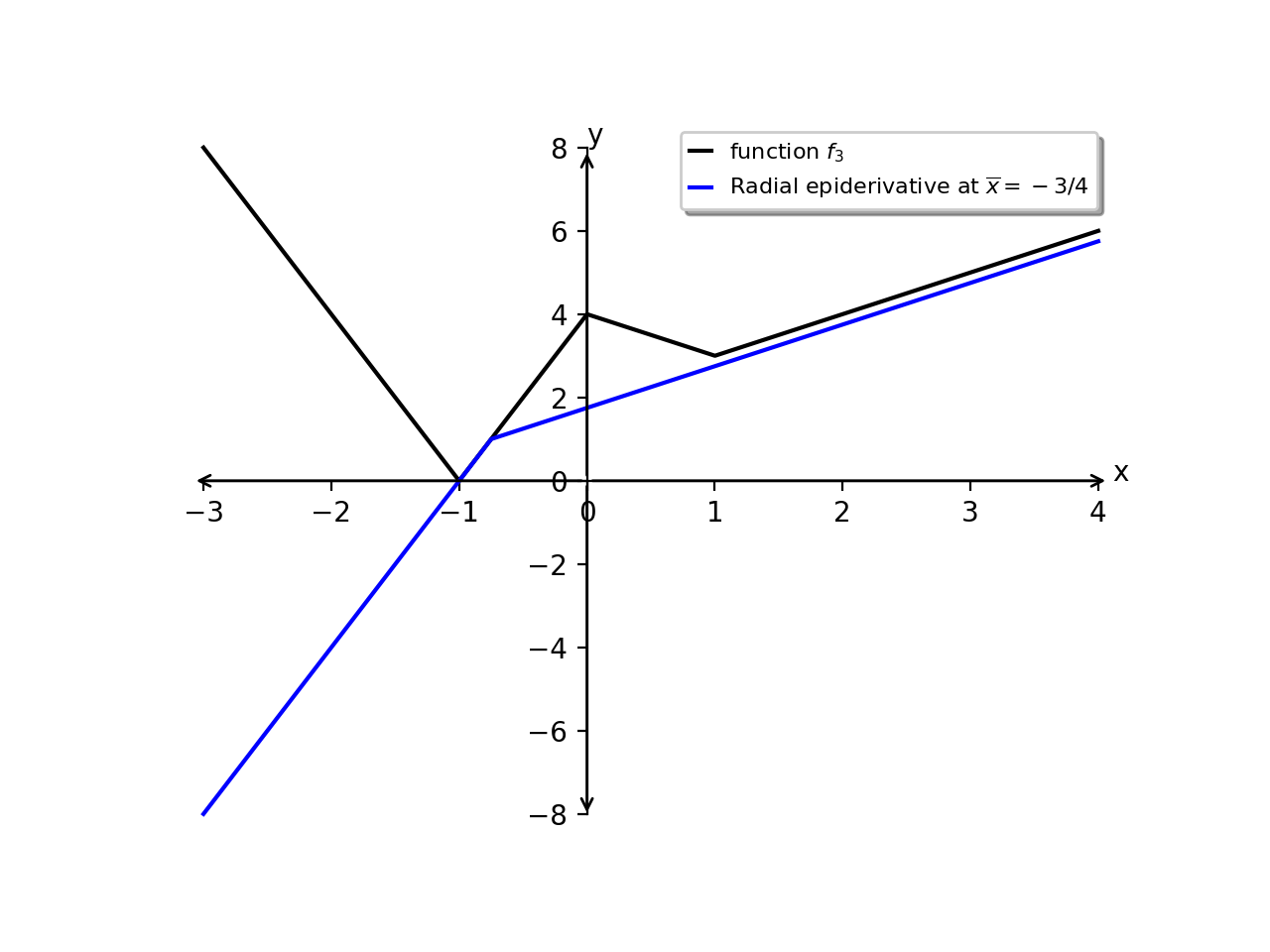}}
		\scalebox{0.39}{\includegraphics{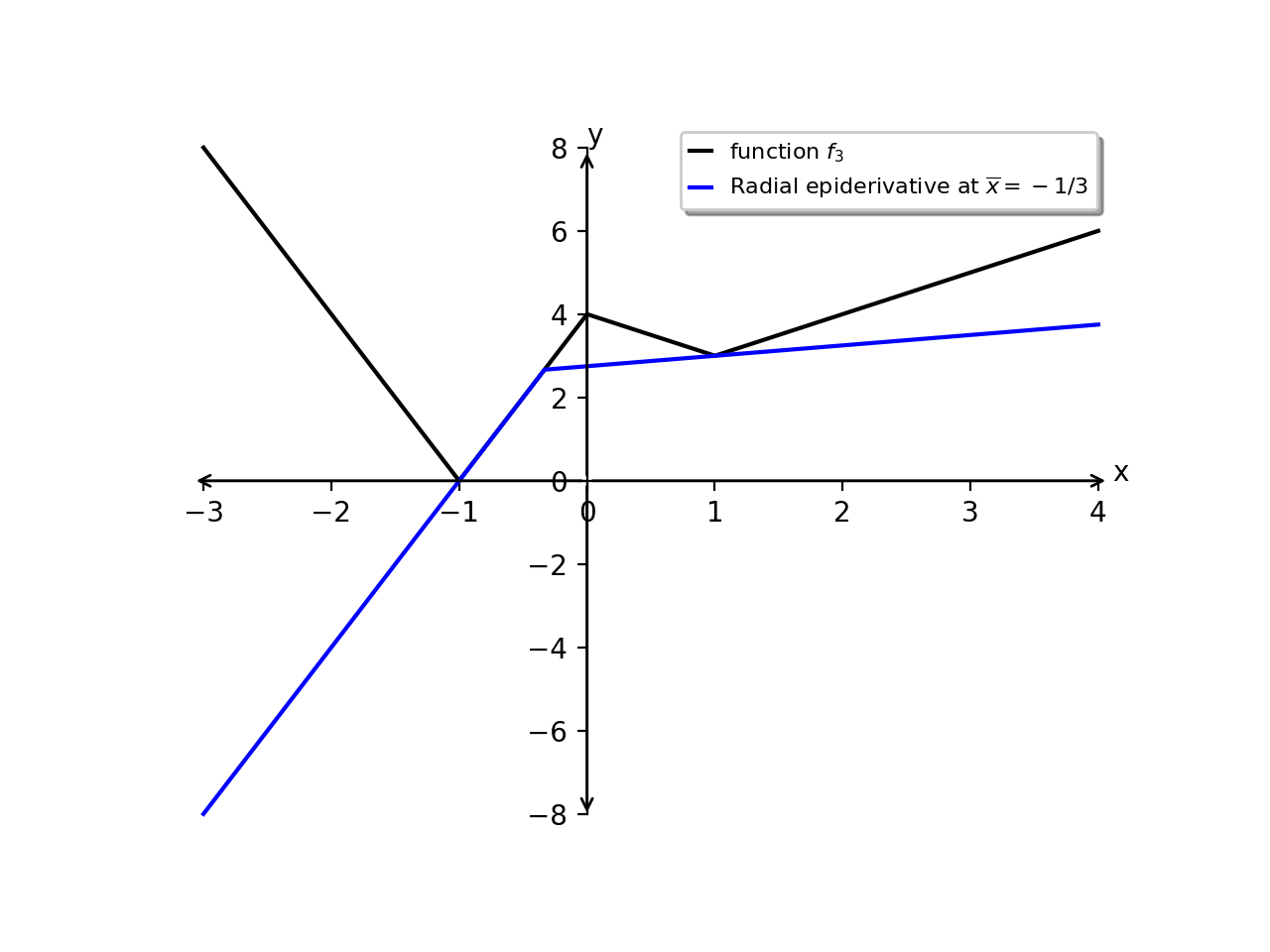}}
		\scalebox{0.39}{\includegraphics{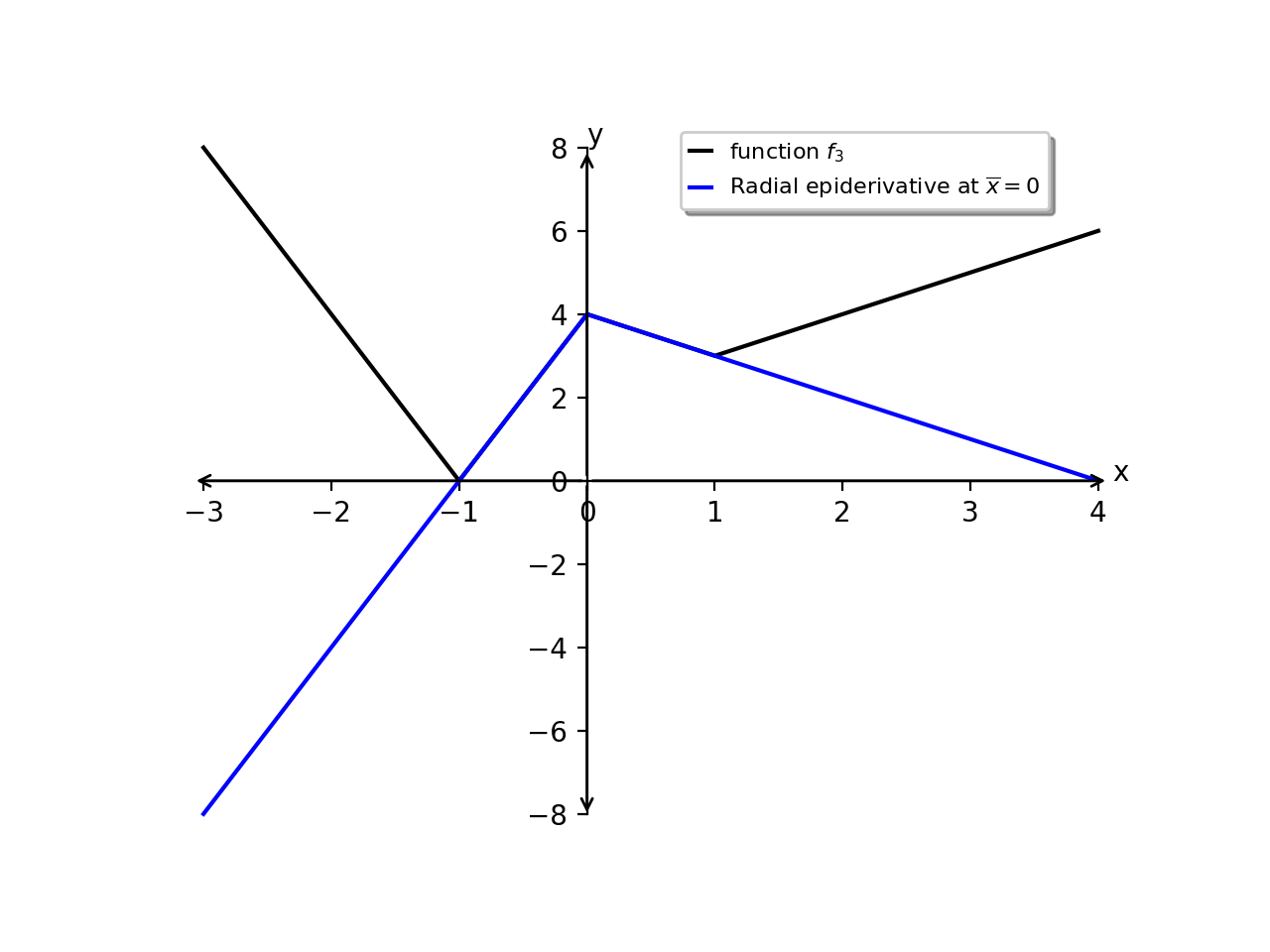}}
		\scalebox{0.39}{\includegraphics{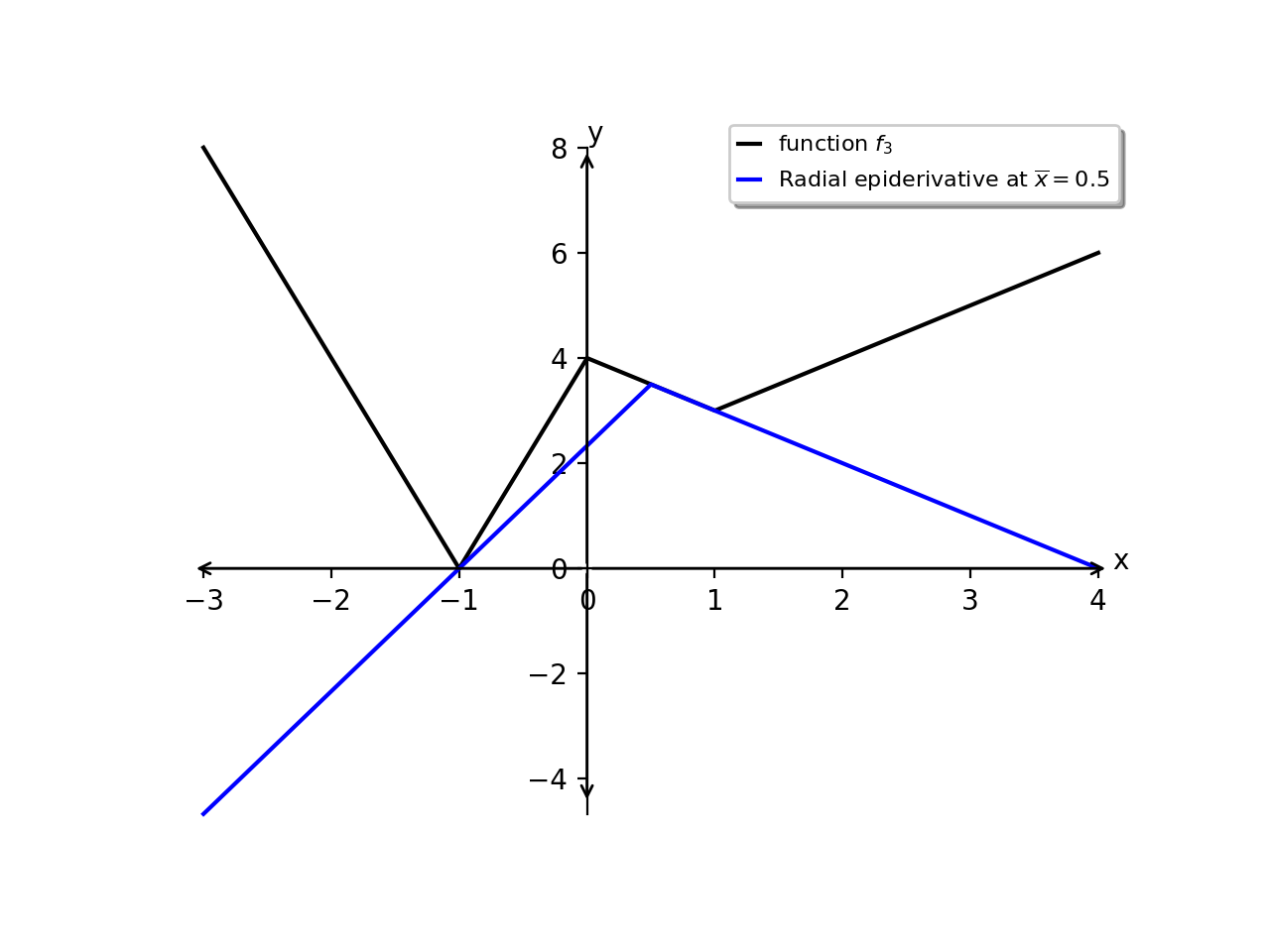}}
		\scalebox{0.39}{\includegraphics{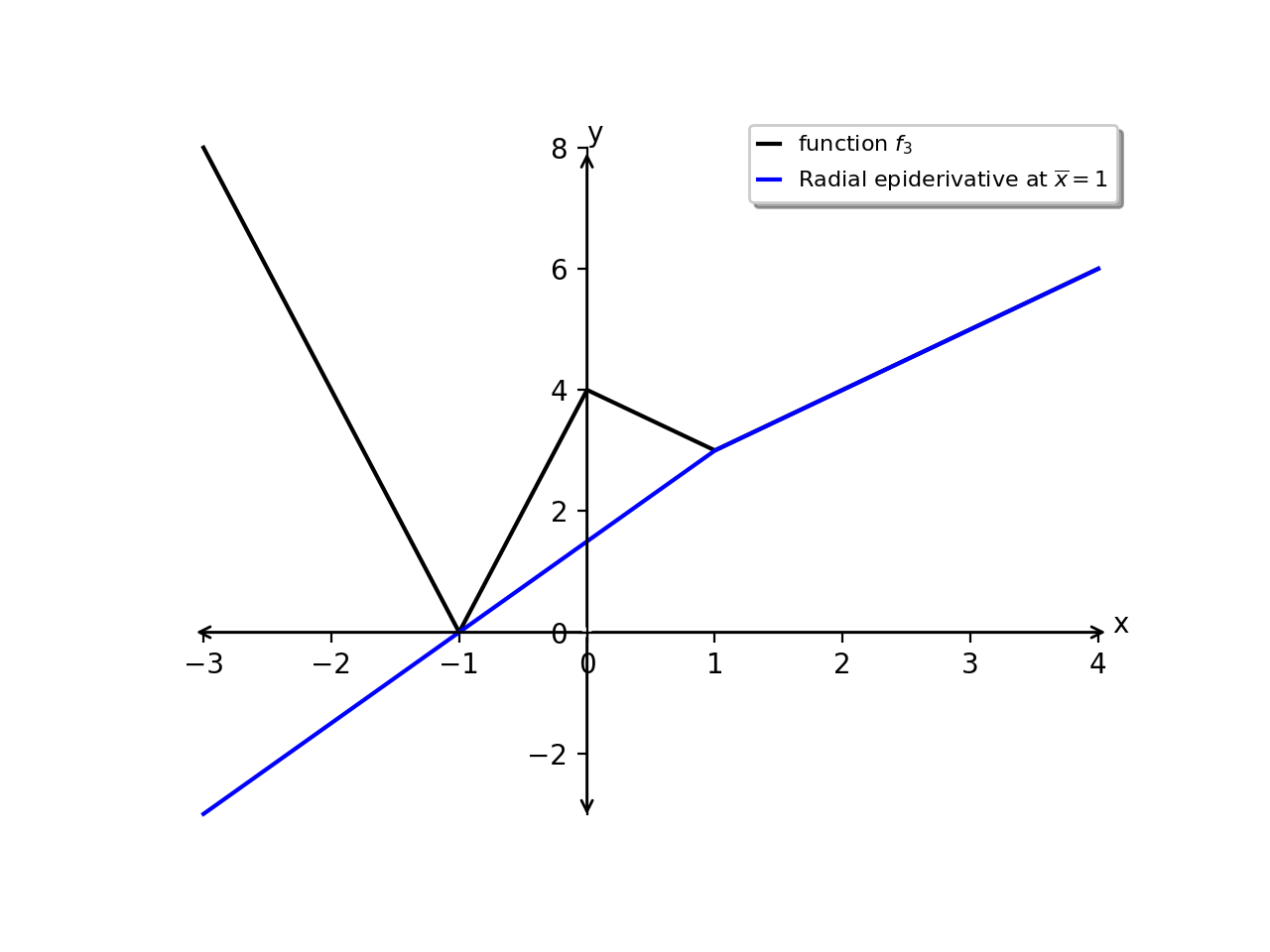}}
		\scalebox{0.39}{\includegraphics{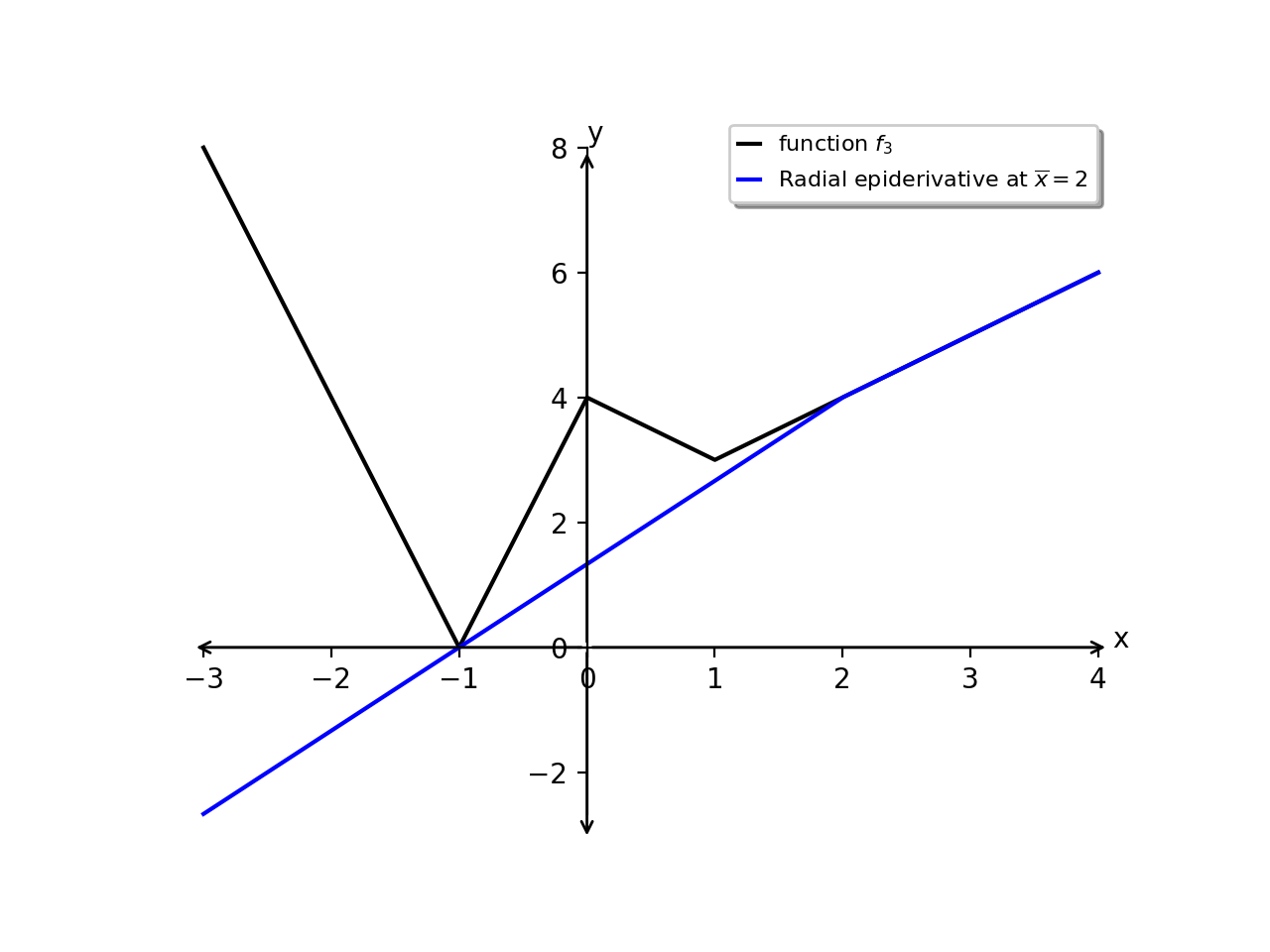}}
		\caption{The graphs of the radial epiderivatives of function $f_3$ at $\overline{x}=-2$, $\overline{x}=-1$, $\overline{x}=-3/4$, $\overline{x}=-1/3$, $\overline{x}=0$, $\overline{x}=1/2$, $\overline{x}=1$, and $\overline{x}=2$ (from left to right)}
		\label{figure:f_3_RE}
	\end{figure}

	\begin{equation}\label{radepiderex3}
	f^r_3(\overline{x};h) =
	\left\{
	\begin{array}{lllllllllll}
	-4h  & \mbox{ if } \overline{x} < -1, h \in \mathbb{R}, \\
	-4h  & \mbox{ if } \overline{x} = -1, h \leq 0, \\
	h    & \mbox{ if } \overline{x} = -1, h > 0, \\
	h    & \mbox{ if } -1 < \overline{x} < -\frac{2}{3}, h>0,\\
	-\frac{(1+4\overline{x})h}{1-\overline{x}}   
	& \mbox{ if } -\frac{2}{3} \leq \overline{x} < 0, h > 0, \\
	4h   & \mbox{ if } -1 < \overline{x} < 0, h<0,\\
	-h   & \mbox{ if } \overline{x} = 0, h > 0, \\
	4h   & \mbox{ if } \overline{x} = 0, h < 0, \\
	\frac{(4-\overline{x})h}{\overline{x}+1}    
	& \mbox{ if } 0 < \overline{x} < 1, h < 0, \\
	-h 	 & \mbox{ if } 0 < \overline{x} < 1, h > 0, \\
	h    & \mbox{ if } \overline{x} = 1, h > 0, \\
	\frac{3h}{2}    
	& \mbox{ if } \overline{x} = 1, h < 0, \\
	h    & \mbox{ if } \overline{x} > 1, h > 0, \\
	\frac{(2+\overline{x})h}{1+\overline{x}} & \mbox{ if } \overline{x} > 1, h < 0. \\
	\end{array}
	\right.
	\end{equation}


\end{example}


\section{Regularity conditions}\label{regularity}

In this section, we investigate the relations between the radial epiderivative, the Clarke's derivative, the Rockafellar's subderivative and the (classical) directional derivative.

The following theorem quoted from \cite{RockafellarW2009}, gives a regularity condition in a nonconvex case and explains a relationship between the directional derivative, the Clarke's directional derivative and the Rockafellar's subderivative.


\begin{theorem} {\normalfont \cite[Theorem 9.16, p. 360]{RockafellarW2009}}\label{subderregular}
	A function $f$ that is finite on an open set $O \subset \mathbb{R}^n$ is both strictly continuous (locally Lipschitz continuous) and regular on $O$ if and only if for every $x \in O$ and $h \in \mathbb{R}^n$ the directional derivative $f^{\prime}(x;h)$ exists, is finite, and depends upper semicontinuously on $x$ for each fixed $h.$ Then $f^{\prime}(x;h)$ depends upper semicontinuously on $x$ and $h$ together and 
	$$
	f^{\prime}(x;h) = \lim_{t\downarrow 0, u \rightarrow h} \frac{f(x+tu) - f(x)}{t} = df(x;h) = f^{\circ}(x;h).
	$$
\end{theorem}

It follows from the definitions of the directional derivative, the Clarke's directional derivative, the subderivative and the radial epiderivative that
\begin{equation}\label{raddirclarkederdir}
f^r(\overline{x};x-\overline{x}) \leq \text{d}f(\overline{x};h) \leq f^{\prime} (\overline{x};x-\overline{x}) \leq f^{\circ} (\overline{x};x-\overline{x})
\end{equation}
for all $x.$

The following theorem provides a condition (different from that given in Theorem \ref{subderregular}), under which a given nonconvex function becomes regular and equality holds in \eqref{raddirclarkederdir}.


\begin{theorem} \label{radepidereqclarkeder}
	Let $(\mathbb{X}, \|\cdot\|_\mathbb{X})$ be a real normed space, let $\mathbb{S}$ be a nonempty subset of the real normed space and let $f:\mathbb{S} \rightarrow \mathbb{R}$ be a
	proper function. Assume that $f$ has a finite directional derivative, a finite Clarke directional derivative and a finite subderivative at 
	$\overline{x} \in \mathbb{X}$ in every direction $x-\overline{x}$ with arbitrary $x
	\in \mathbb{X}.$ Then $f$ is radially epidifferentiable at	$\overline{x}$ and
	\begin{equation}\label{radekclarkederdir}
	f^r(\overline{x};x-\overline{x})= \text{d}f(\overline{x};x-\overline{x}) = f^{\prime} (\overline{x};x-\overline{x})  = f^{\circ} (\overline{x};x-\overline{x})
	\end{equation}
	if and only if
	\begin{equation}\label{maindircondeqthm}
	f(x)-f(\overline{x})\geq  f^{\circ} (\overline{x};x-\overline{x})
	\quad \mbox{ for all } x \in \mathbb{X}.
	\end{equation}
\end{theorem}	


\begin{proof}
	Assume that $f$ has a finite Clarke directional derivative at
	$\overline{x} \in \mathbb{X}$ in every direction $x-\overline{x}$ with arbitrary $x
	\in \mathbb{X}$ and that \eqref{maindircondeqthm} is satisfied. 
	Since $f^{\circ}(\overline{x};\cdot)$ is a Lipschitz function by \cite[Proposition 2.1.1, p. 25]{Clarke1983}, it is also lower Lipschitz. Then, it follows from \eqref{maindircondeqthm} that $f$ is lower Lipschitz at $\overline{x} \in \mathbb{X},$ too. Hence by Theorem \ref{radepiderlowerlip}, $f$ is radially epidifferentiable at
	$\overline{x} \in \mathbb{X}$ in every direction $x-\overline{x}$ with arbitrary $x
	\in \mathbb{X}.$ 	
	By using the representation \eqref{radepilimdef} for the radial epiderivative, we have:
	\begin{equation*} \label{radepilimdef2}
	f^{r}(\overline{x};h) = \inf_{t > 0} \lim_{ u \rightarrow h} \frac{f(\overline{x}+ tu)- f(\overline{x})}{t}\geq \inf_{t > 0} \lim_{ u \rightarrow h} \frac{f^{\circ}(\overline{x}; tu)}{t} =f^{\circ}(\overline{x}; h)
	\end{equation*}
	for all $h \in \mathbb{X},$ where the last equality above, is obtained due to the positive homogeneity and the Lipschitz continuity of $f^{\circ}(\overline{x}; \cdot)$ \cite[Proposition 2.1.1 (a),(b), p. 25]{Clarke1983}. Then, the claim follows from \eqref{raddirclarkederdir}.
	
	Now assume that \eqref{radekclarkederdir} is satisfied for all $x \in \mathbb{X}.$ Then the claim follows from the inequality 
	\begin{equation*}
	f(x)-f(\overline{x}) \geq f^{r}(\overline{x};x-\overline{x})
	\end{equation*}
	for all $x \in \mathbb{X},$ which is satisfied for the radial epiderivative due to its definition  (see \eqref{radial20}). $\Box$
\end{proof}


It follows from the definitions of the generalized derivatives mentioned in Theorem \ref{radepidereqclarkeder}, there may be cases when the directional derivative and the subderivative does exist but the Clarke derivative does not exist. Moreover there may be a case when the subderivative does exist but the directional derivative does not exist. The following corollaries which are straightforward from Theorem \ref{radepidereqclarkeder}, provide equality relations between the existing derivatives.


\begin{corollary} \label{radepidereqdirder}
	Let $(\mathbb{X}, \|.\|_\mathbb{X})$ be real normed space and let $f:\mathbb{R} \rightarrow \mathbb{R}$ be a
	proper function. Assume that $f$ has both the finite directional derivative and the finite subderivative at
	$\overline{x} \in \mathbb{X}$ in every direction $x-\overline{x}$ with arbitrary $x
	\in \mathbb{X}.$ Then $f$ is radially epidifferentiable at	$\overline{x}$ and 
	\begin{equation}\label{radekdirder}
	f^r(\overline{x};x-\overline{x})= \text{d}f(\overline{x};x-\overline{x}) = f^{\prime} (\overline{x};x-\overline{x})
	\end{equation}
	if and only if
	\begin{equation}\label{condradepidereqdirder}
	f(x)-f(\overline{x})\geq  f^{\prime} (\overline{x};x-\overline{x})
	\quad \mbox{ for all } x \in \mathbb{X}.
	\end{equation}
\end{corollary}	


\begin{proof} By the hypothesis, we have:
	$$
	\inf_{t > 0} \lim_{ u \rightarrow h} \frac{f(\overline{x}+ tu)- f(\overline{x})}{t}\geq \inf_{t > 0} \lim_{ u \rightarrow h} \frac{f^{\prime}(\overline{x}; tu)}{t} =f^{\prime}(\overline{x}; h).
	$$
	This means by the definition \eqref{radepilimdef} of the radial epiderivative that $f$ is radially epidifferentiable at	$\overline{x}.$ Then, \eqref{radekdirder} follows from \eqref{raddirclarkederdir}.
	The proof of the second part is similar to that of Theorem \ref{radepidereqclarkeder}. $\Box$
\end{proof}


\begin{corollary} \label{radepidereqsubder}
	Let $(\mathbb{X}, \|.\|_\mathbb{X})$ be real normed space, et $\mathbb{S}$ be a nonempty subset of the real normed space and let $f:\mathbb{S} \rightarrow \mathbb{R}$ be a
	proper function. Assume that $f$ has the finite subderivative at
	$\overline{x} \in \mathbb{X}$ in every direction $x-\overline{x}$ with arbitrary $x
	\in \mathbb{X}.$ Then  $f$ is radially epidifferentiable at	$\overline{x}$ and
	\begin{equation*}\label{radeksubder}
	f^r(\overline{x};x-\overline{x})= \text{d}f(\overline{x};x-\overline{x})
	\end{equation*}
	if and only if
	\begin{equation}\label{condradepidereqsubder}
	f(x)-f(\overline{x})\geq  \text{d}f(\overline{x};x-\overline{x})
	\quad \mbox{ for all } x \in \mathbb{X}.
	\end{equation}
\end{corollary}	


\begin{proof} The proof is similar to the proof of Corollary \ref{radepidereqdirder}. $\Box$
\end{proof}


Now consider illustrative examples for the above assertions.


\begin{example}
	Consider the function 
	\begin{equation*}
	f_1(x) =
	\left\{
	\begin{array}{ll}
	-x+3    & \mbox{ if } x < 1, \\
	x & \mbox{ if } x \geq 1
	\end{array}
	\right.
	\end{equation*}
	from Example \ref{example1}, where it has been shown that $f^{r}_1(1;h)=|h|.$  It is easy to see that 
	\begin{equation*}\label{dirder1ex1}
	\text{d}f_1(1;h) =f^{\prime}_1(1;h) =f^{\circ}_3(1;h)=
	\left\{
	\begin{array}{ll}
	+\infty  & \mbox{ if } h < 0, \\
	h & \mbox{ if } h>0.
	\end{array}
	\right.
	\end{equation*}
	Note that this example nicely illustrates the Theorems \ref{subderregular} and \ref{radepidereqclarkeder} for different choices of $\overline{x}$ and $h.$
	$\Box$
\end{example}


The following example demonstrates that the subderivative may exist, when both the directional derivative and the Clarke directional derivative do not exist. This example again demonstrates the case when the subderivative and the radial epiderivative are equal.

\begin{figure}[h]
	\centering
	\scalebox{0.4}{\includegraphics{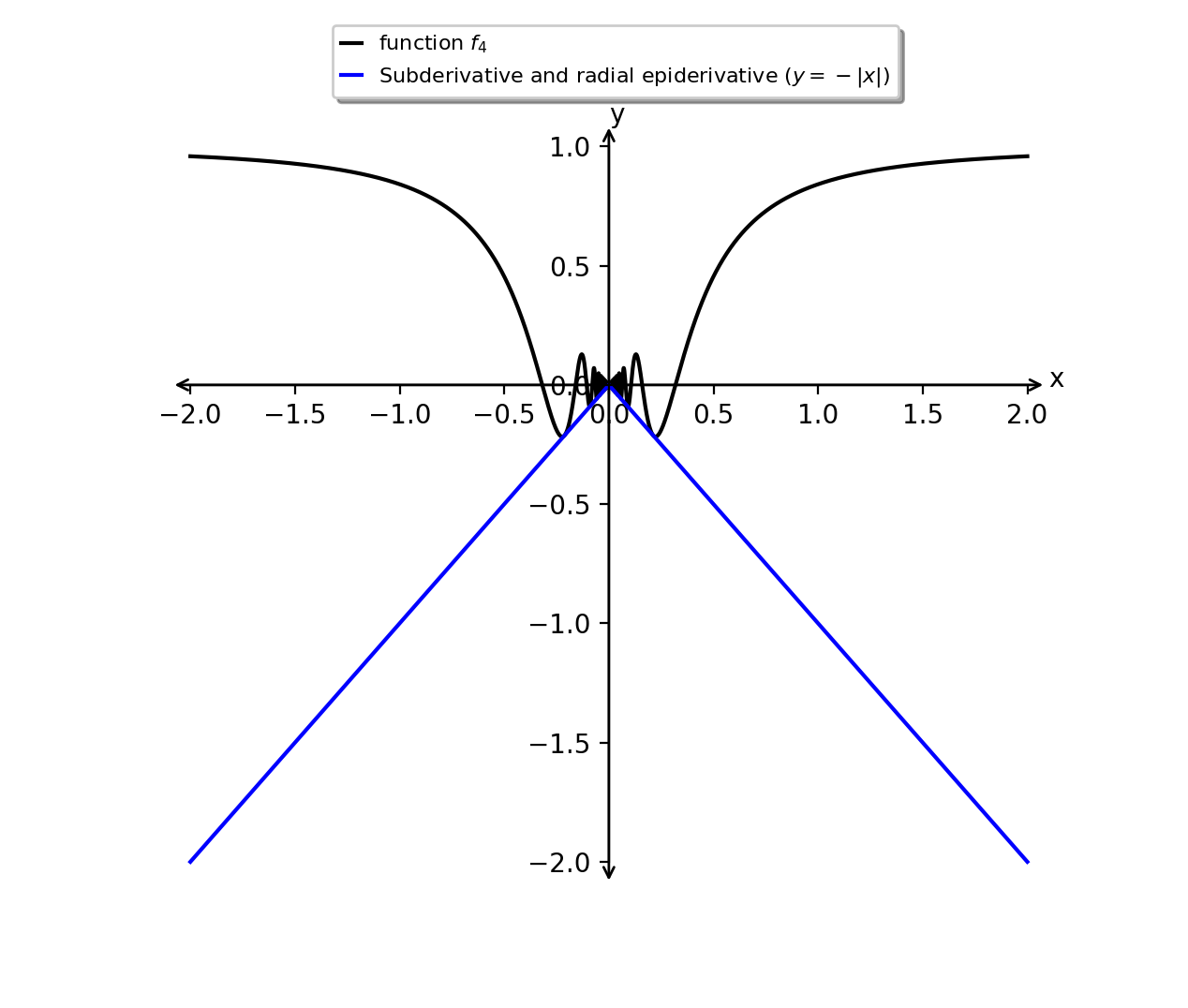}}
	\caption{The graph of the radial epiderivative and the subderivative  of $f_4$ at $\overline{x}=0.$}
	\label{figure:f_4}
\end{figure}

\begin{example}\label{allder1}
	Let  
	\begin{equation*}\label{funcex4}
	f_4(x) =
	\left\{
	\begin{array}{ll}
	x\sin (\frac{1}{x})  & \mbox{ if } x \neq 0, \\
	0                    & \mbox{ if } x =0
	\end{array}
	\right.
	\end{equation*}
	whose graph is depicted in Figure \ref{figure:f_4}. Then the conditions of Corollary \ref{radepidereqsubder} are satisfied and
	
	\begin{equation*}\label{radepiderex4}
	f^r_4 (0;h) = \text{d}f_4 (0;h)= \left\{
	\begin{array}{ll}
	h  & \mbox{ if } h < 0, \\
	-h  & \mbox{ if } h > 0.
	\end{array}
	\right.
	\end{equation*}
It is easy to show that both the directional derivative $f^{\prime}_4 (0;h)$ and the Clarke	directional derivative $f^{\circ}_{5} (0;h)$ fail to exist.

\end{example}

\begin{example}\label{xsinlnx}
	Consider the following function (see also \cite[Fig 6-4, p.199]{RockafellarW2009}):
	\begin{equation*}\label{funcexxsinlnx}
	f_{5}(x) =
	\left\{
	\begin{array}{ll}
	x \sin (\ln |x|)  & \mbox{ if } x \neq 0, \\
	0                    & \mbox{ if } x =0
	\end{array}
	\right.
	\end{equation*}
Then
	
	\begin{equation*}\label{radepiderexxsinlnx}
	f^r_{5} (0;h) = \left\{
	\begin{array}{ll}
	h  & \mbox{ if } h < 0, \\
	-h  & \mbox{ if } h > 0.
	\end{array}
	\right.
	\end{equation*}
	
	Note that the conditions of Corollary \ref{radepidereqsubder} are satisfied and hence $f^r_{5} (0;h) = \text{d}f_{5} (0;h)=-|h|.$ On the other hand $f^{\prime}_{5} (0;h) $ does not exist, but  $f^{\circ}_{5} (0;h) = |h|$ (See Figure \ref{figure:f_5}).
	
	\begin{figure}[h]
		\scalebox{0.39}{\includegraphics{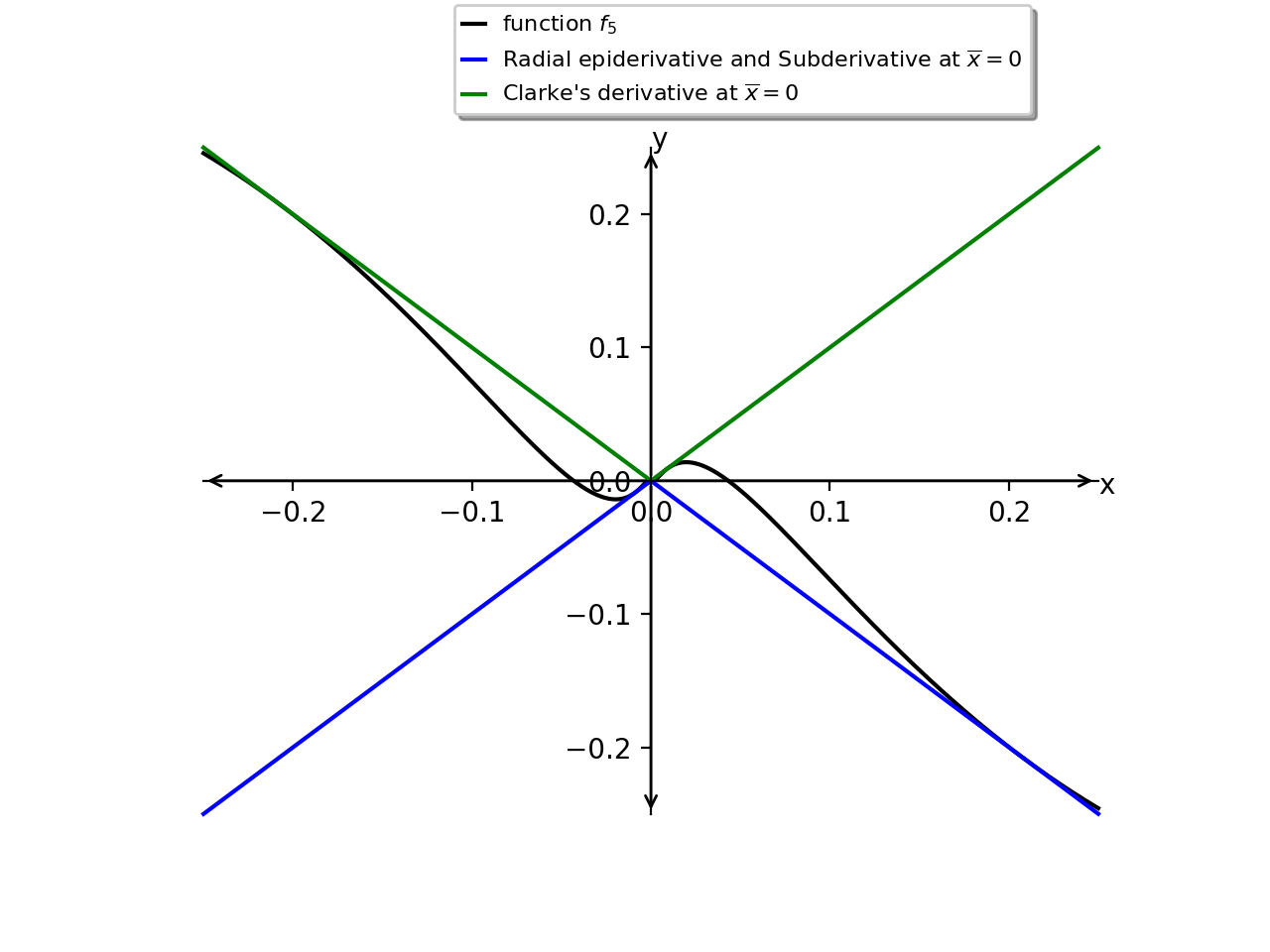}}
		\scalebox{0.39}{\includegraphics{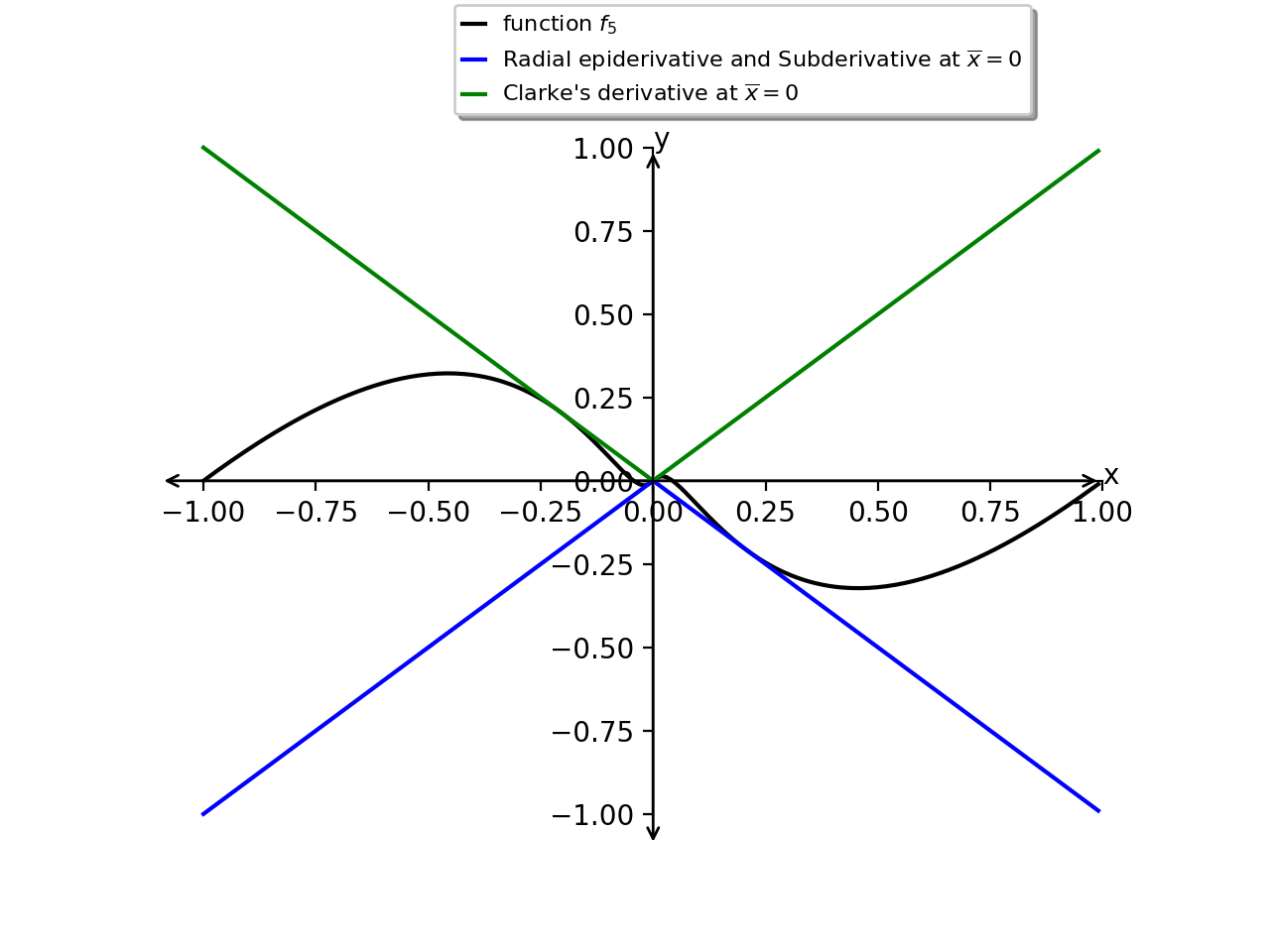}}
		\caption{The graph of radial epiderivative, subderivative and Clarke's derivative of function $f_5$ at $\overline{x}=0$ in the interval $x \in [-0.25,0.25] $ (left) and $x \in [-1,1] $ (right)}
		\label{figure:f_5}
	\end{figure}

\end{example}


\begin{example}\label{x2sin1overx}
	Consider the following function (see also \cite[Execcise 8.8, p.304]{RockafellarW2009} and \cite[Example 2.2.3, p.33]{Clarke1983}):
	\begin{equation*}\label{funcx2sin1overx}
	f_{6}(x) =
	\left\{
	\begin{array}{ll}
	x^2\sin (\frac{1}{x})  & \mbox{ if } x \neq 0, \\
	0                    & \mbox{ if } x =0
	\end{array}
	\right.
	\end{equation*}
	
	For $\overline{x}=0$ we have $f_{6}^{\prime}(0;h) = 0$ for all $h \in \mathbb{R}. $ Since the derivative mapping $\nabla f_{6}$ is discontinuous at $\overline{x} =0,$ the conditions of Theorem \ref{subderregular}, are not satisfied. In a similar way, for the subderivative we obtain that $\text{d}f_{6} (0;h) = 0.$ On the other hand, the condition \eqref{mainconddirder} of Corollary \ref{radepidereqdirder} is also not satisfied. For the radial epiderivative at $\overline{x} =0,$  we have:
	\begin{equation*}\label{radfuncx2sin1overx}
	f^r_{6}(0;h) =
	\left\{
	\begin{array}{ll}
	h              & \mbox{ if } h < 0, \\
	-\alpha h      & \mbox{ if } h \geq 0
	\end{array}
	\right.
	\end{equation*}
	where $-1<-\alpha = x_0^2\sin (\frac{1}{x_0}) <0,$ with $0.2 < x_0 < 0.3$ (see Figure \ref{figure:f_6}).
	
	For the Clarke directional derivative of $f_{6}$ at $\overline{x}=0$ we have $f_{6}^{\circ}(0;h) = |h| $ for all $h$ (see, \cite[Example 2.2.3, p.33]{Clarke1983} ). All derivatives for this example, are depicted in Figure \ref{figure:f_6}. 
	
	\begin{figure}[h]
		\centering
		\scalebox{0.34}{\includegraphics{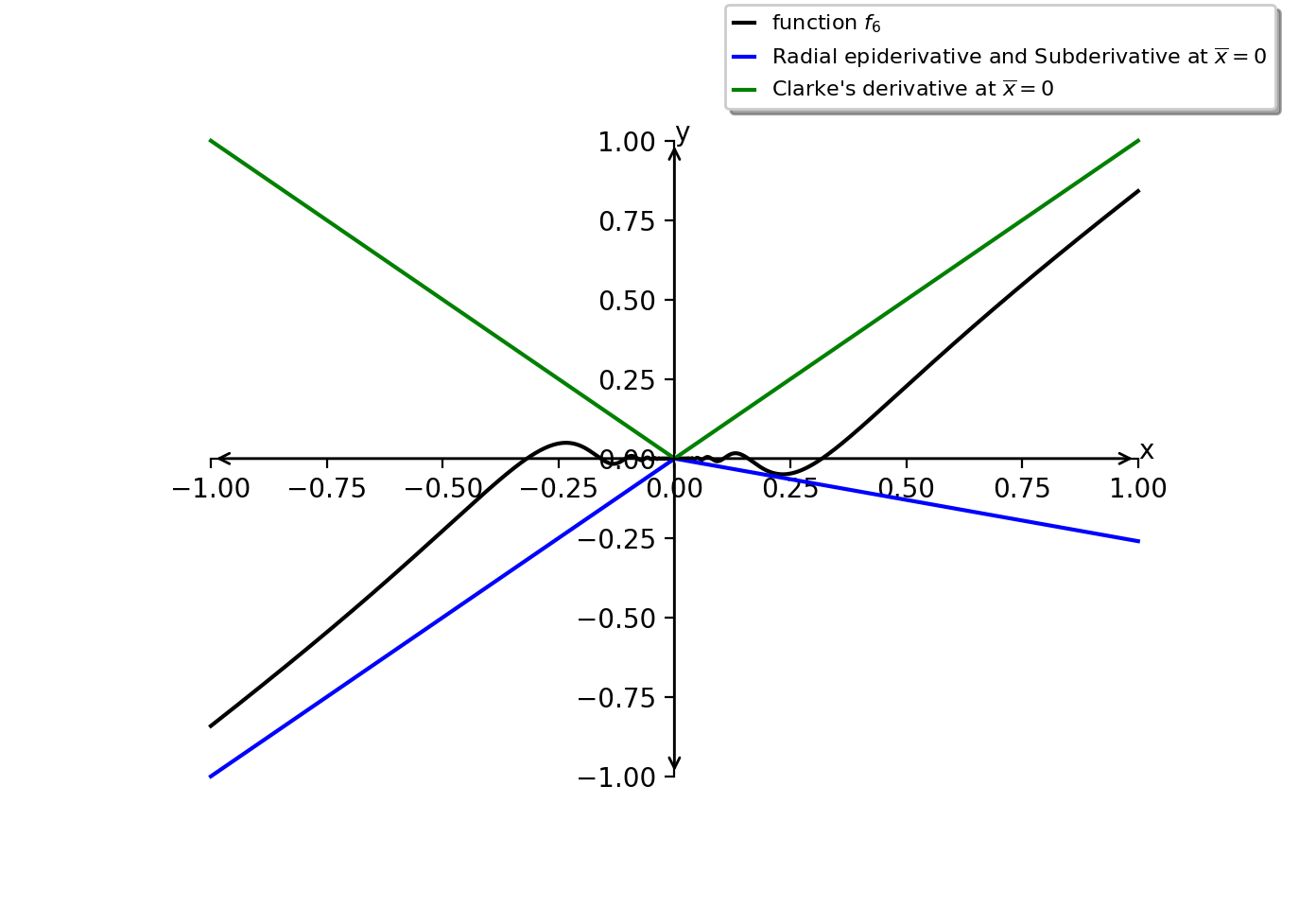}}
		\caption{The graph of radial epiderivative with $-\alpha=-0.26$, subderivative and Clarke's derivative of function $f_6$ at $\overline{x}=0$}
		\label{figure:f_6}
	\end{figure}	
	
	This function demonstrates the case, when the given function has all the generalized derivatives considered in this paper, with
	
	$$
	f^r_{6} (0;h) \neq \text{d}f_{6} (0;h) = f^{\prime}_{6}(0;h)= 0 \neq f^{\circ}_{6} (0;h) = |h|.
	$$
	
\end{example}


\begin{example}\label{x2sin21overx}
	Let  
	\begin{equation*}\label{funcx2sin21overx}
	f_{7}(x) =
	\left\{
	\begin{array}{ll}
	x^2\sin^2 (\frac{1}{x})  & \mbox{ if } x \neq 0, \\
	0                        & \mbox{ if } x =0
	\end{array}
	\right.
	\end{equation*}
	whose graph is depicted in Figure \ref{figure:f_7}. 
	\begin{figure}[h]
		\centering
		\scalebox{0.4}{\includegraphics{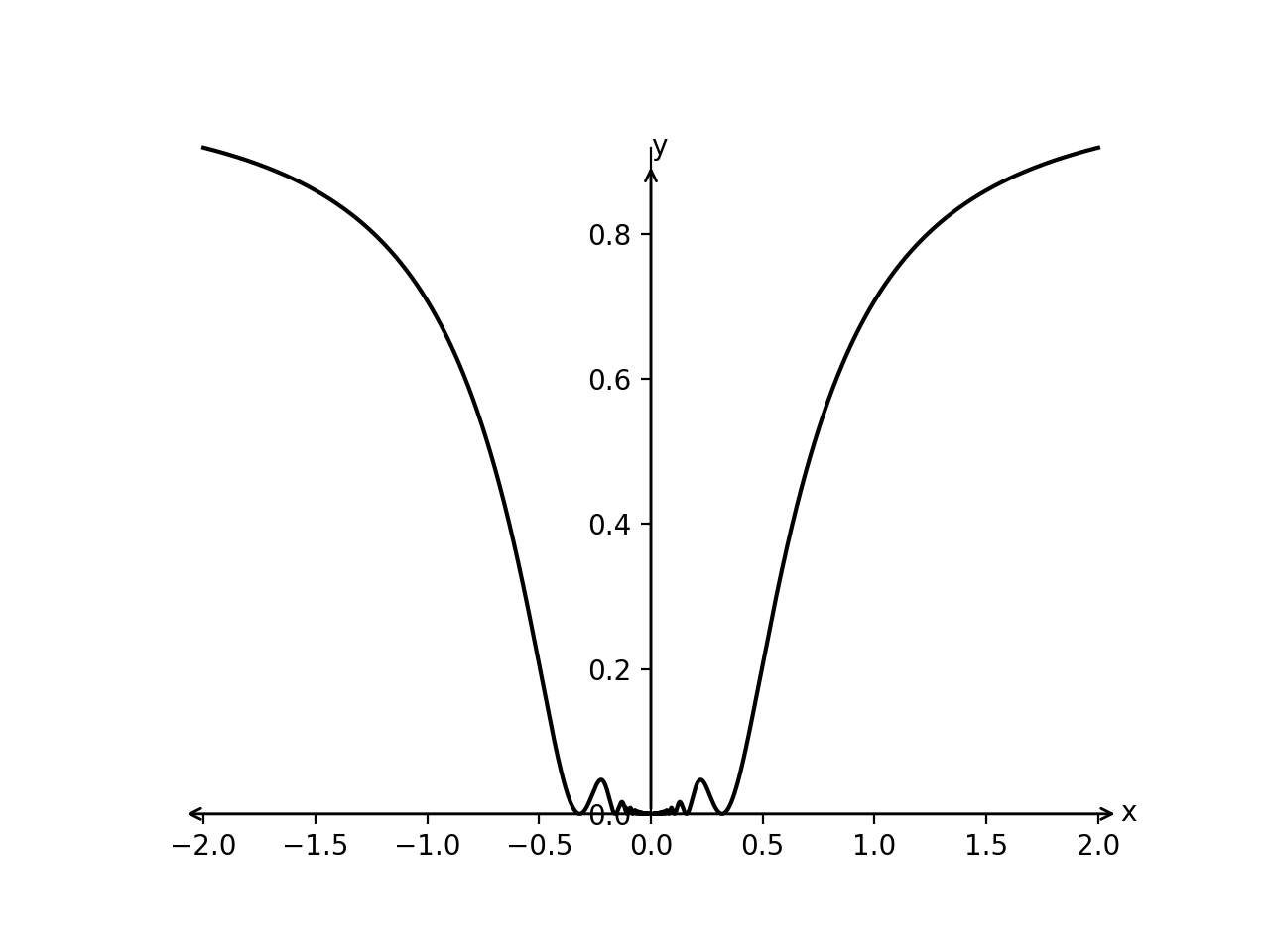}}
		\caption{The graph of function $f_7$}
		\label{figure:f_7}
	\end{figure}	
\end{example}
For $\overline{x}=0$ we have $f_{7}^{\prime}(0;h) =0$
for all $h \in \mathbb{R}, $ but the derivative mapping $\nabla f_{7}$ is discontinuous at $\overline{x} =0$ and the conditions of Theorem \ref{subderregular}, are not satisfied. On the other hand, we have $f_{7}(x) - f_{7}(0) \geq f_{7}^{\prime}(0;x-0)$ for all $x \in \mathbb{R},$ which shows that the conditions of Corollary \ref{radepidereqdirder} are satisfied and as a result, we have:
$$
f^r_{7} (0;h) = \text{d}f_{7} (0;h) = f^{\prime}_{7}(0;h)=0. 
$$ 

This function also demonstrates the case when the given function is directionally differentiable, but is not Clarke directionally differentiable.


\begin{example}\label{allderf3func}
	Consider the function $f_3$ from Example \ref{example3} (see Figure \ref{figure:f_3_RE}) and try to interpret Theorem \ref{radepidereqclarkeder} and Corollary \ref{radepidereqdirder}. Recall that
	
	\begin{equation*}
	f_3(x) =
	\left\{
	\begin{array}{ll}
	4|x+1|    & \mbox{ if } x \leq 0, \\
	|x-1| + 3 & \mbox{ if } x > 0.
	\end{array}
	\right.
	\end{equation*}
	Then, all conditions of Theorem \ref{radepidereqclarkeder} and the assumption \eqref{maindircondeqthm} are satisfied at points $\overline{x} < -1$ and
	\begin{equation*}\label{radekclarkedirderrock}
	f_3^r(\overline{x};h)= \text{d}f_3(\overline{x};h) = f_3^{\prime} (\overline{x};h) = f_3^{\circ} (\overline{x};h)
	\end{equation*}
	for all $\overline{x} < -1$ and  $h	\in \mathbb{R}.$
	For $\overline{x} = -1$ we have:
	\begin{equation*}\label{clarkedirderrock}
	f_3^r(-1;h)\neq \text{d}f(-1;h) = f_3^{\prime} (-1;h) = f_3^{\circ} (-1;h)=4h
	\end{equation*}
	for $h>0,$ where the assumption \eqref{maindircondeqthm} is not satisfied at $\overline{x} = -1.$
	Finally, it is clear that the condition \eqref{condradepidereqdirder} of Corollary \ref{radepidereqdirder} is satisfied at the point $\overline{x} = 0,$ where the condition of Theorem \ref{subderregular} is not satisfied and as a result  we have:
	\begin{equation*}\label{raddirderrock}
	f_3^r(0;h)= \text{d}f_3(0;h) = f_3^{\prime} (0;h)=\left\{
	\begin{array}{ll}
	-4h   & \mbox{ if } h \leq 0, \\
	h & \mbox{ if } h > 0.
	\end{array}
	\right. \neq f_3^{\circ}(0;h)
	\end{equation*}
	for all $h	\in \mathbb{R},$ where
	\begin{equation*}
	f_3^{\circ}(0;h) =
	\left\{
	\begin{array}{ll}
	h   & \mbox{ if } h \leq 0, \\
	4h & \mbox{ if } h > 0.
	\end{array}
	\right.
	\end{equation*}

	$\Box$
\end{example}


\begin{remark}\label{tractabilityofrdepider}
	An extension of the optimaity condition \eqref{pshrockoptcond} to a nonconvex case established in terms of the weak subgradients \cite[Theorem 4]{KasimbeyliM2011}, the representation of the radial epiderivative as a support function of the  weak subdifferential set in nonconvex case (see \cite[Theorems 4.5 and 4.6]{KasimbeyliM2009}), Theorem \ref{radepiderlowerlip} on the necessary and sufficient condition for the radial epidifferentiability, regularity relations given in Theorem  \ref{radepidereqclarkeder}, Corollaries \ref{radepidereqdirder} and \ref{radepidereqsubder} as well as the illustrative examples presented above help to better understand the radial epiderivative notion. On the other hand, Theorem  \ref{radepidereqclarkeder}, Corollaries \ref{radepidereqdirder} and \ref{radepidereqsubder} show that computational methods and approaches existing for computing the directional derivative, subderivative and the generalized derivative can be used to estimate the radial epiderivative.
	\end{remark}


\begin{remark}\label{generalizationofsupportrel}
The regularity conditions \eqref{maindircondeqthm}, \eqref{condradepidereqdirder} and \eqref{condradepidereqsubder} given in Theorem \ref{radepidereqclarkeder}, Corollaries \ref{radepidereqdirder} and \ref{radepidereqsubder} respectively, not only provide necessary and sufficient conditions for equality of the generalized derivatives under consideration, but these conditions also generalize the main affine support relation \eqref{gradineq} of the convex analysis to a nonconvex case. The examples considered above, explain and illustrate that these conditions are valid not only for convex but also for nonconvex functions.
\end{remark}


\section{Computing weak subgradients by using radial epiderivatives}\label{computingweaksubgr}

This section presents two theorems with constructive proofs, which provide explicit formulas for approximate computing weak subgradients at a given point. These theorems can be considered as versions of the same theorem for the Euclidean ($\ell_2$)  and the $\ell_1-$ norms.

We first give the following lemma which  plays an important role in the proof of the subsequent theorems.


\begin{lemma}\label{weaksubdifrad}
	Let $f:\mathbb{R}^n\rightarrow \mathbb{R}$ be a proper  function, $\overline{x} \in \mathbb{R}^n$, and $f(\overline{x})$ be finite. If $f$ has a radial epiderivative $f^{r}(\overline{x};x)$ for every $x \in \mathbb{R}^n,$
	then $f^{r}(\overline{x};\cdot)$ is weakly subdifferentiable at $0_{\mathbb{R}^n,}$  $f$ is weakly subdifferentiable at $\overline{x}$ and
	\begin{equation}\label{partialsubset}
	\partial^{w}f^{r}(\overline{x};0) = \partial^{w} f(\overline{x}).
	\end{equation}
\end{lemma}
\unskip

\begin{proof} Suppose that the radial epiderivative $f^{r}(\overline{x};\cdot)$ exists 	and is given by 
	\eqref{radial20}. It follows from this relation that $f^{r}(\overline{x};\cdot)$ is bounded from below on some neighborhood of $0_{\mathbb{R}^n}$ and that, it is
	a positively homogeneous function (see Lemma \ref{radposhom}). Therefore by Theorem \ref{wsubdiffunc}, it is weakly subdifferentiable at $0_{\mathbb{R}^n}$.
	
	On the other hand, it follows from the definition of the radial cone that
	$$
	(x-\overline{x},f(x)-f(\overline{x}))\in \mathbb{R}(\mathit{epi}(f);(\overline{x},f\overline{x})).
	$$
	Indeed, the radial cone $\mathbb{R}(\mathit{epi}(f);(\overline{x},f\overline{x}))$ consists of elements of the form 
	$$\lim_{n \rightarrow +\infty}t_n[(x_n,y_n) - (\overline{x},f(\overline{x}))],$$ 
	where in particular we can take $t_n =1, x_n =x, y_n = f(x)$ for all $n,$ which leads to the element $(x-\overline{x},f(x)-f(\overline{x})).$ Therefore  by (\ref{radial20}), we obtain the following relation:
	\begin{equation}\label{mainradcond}
	f^{r}(\overline{x};x-\overline{x})\leq f(x)-f(\overline{x})
	\end{equation}
	for all $x \in \mathbb{R}^n.$
	
	Now we show that $f$ is weakly subdifferentiable at $\overline{x}$ and that \eqref{partialsubset} is satisfied.
	Let $(x^*,c)\in \partial^{w} f^{r}(\overline{x};0).$ Then
	\[
	f^{r}(\overline{x};h) \geq \langle x^*,h \rangle - c\|h\| \quad \mbox{for all } h\in \mathbb{R}^n,
	\]
	or
	\[
	f^{r}(\overline{x};x-\overline{x}) \geq \langle x^*,x-\overline{x} \rangle - c\|x-\overline{x}\| \quad \mbox{for all } x\in \mathbb{R}^n.
	\]
	By \eqref{mainradcond} this implies
	\[
	f(x)-f(\overline{x}) \geq \langle x^*,x-\overline{x} \rangle - c\|x-\overline{x}\|, \quad \mbox{for all } x\in \mathbb{R}^n,
	\]
	which means that $(x^*,c)\in \partial^{w} f(\overline {x}).$

	If $(x^*,c) \in \partial^{w}f(\overline{x}),$ then for any
	fixed $h \in \mathbb{R}^n$ we have:
	\[
	f^{r} (\overline{x}; h) =
	\inf_{t > 0} \lim_{ u \rightarrow h} \frac{f(\overline{x}+ tu)- f(\overline{x})}{t}
	\]
	\vspace{-12pt}
	\[
	\geq \inf_{t > 0} \lim_{ u \rightarrow h} \frac{\langle x^{*},tu \rangle - c\|
		tu\|}{t} = \langle x^{*},h \rangle - c\| h \|
	\]
	that is $(x^*,c) \in \partial^{w} f^{r}(\overline{x};0)$
	and hence the proof is completed. $\Box$
\end{proof}


\begin{theorem} \label{computingvall2}
	Let $f:\mathbb{R}^n\rightarrow \mathbb{R}$ be a proper  function, $\overline{x} \in (\mathbb{R}^n,\|\cdot\|_2)$,
	and let $\overline{y}=f(\overline{x})$ be finite. Assume that, $f$ is radially epidifferentiable at $\overline{x}.$ Then, for every $\varepsilon > 0$ and $h \in \mathbb{R}^n \setminus \{0_{\mathbb{R}^n}\}$ with $\|h\| =1,$ there exists a weak subgradient  $(v,c)  \in \partial^wf(\overline{x})$ such that
	\begin{equation}\label{weaksubdircomp}
	v=(c+f^r(\overline{x};h)- \varepsilon)h
	\end{equation}
	or 
	\begin{equation}\label{radialepidercomp}
	f^r(\overline{x};h) = \langle v,h \rangle - c + \varepsilon.
	\end{equation}
	
\end{theorem}


\begin{proof}
	Let $h \in \mathbb{R}^n$ be an arbitrary point with $\|h\| =1$. By the positive homogeneity of $f^r(\overline{x};\cdot)$  (see Lemma \ref{radposhom}), it is
	sufficient to consider elements $h$ with $\|h\|=1$. Let $\varepsilon > 0$ be an arbitrary positive number. We will show that, there exist a nonnegative number $c$, and a vector
	$v \in \mathbb{R}^n$ (possibly depending on $c$ and $h$) such that the pair $(v,c)$ is a weak subgradient of $f^r(\overline{x};\cdot)$ at zero, that is the following inequality is satisfied for every $x \in \mathbb{R}^n:$
	
	\begin{equation*}
	f^r(\overline{x};x) - f^r(\overline{x};0) \geq  \langle v,x-0 \rangle -c\left\| x-0\right\|.
	\end{equation*}
	Since $f^r(\overline{x};0)=0,$ the above inequality can be written simply in the following form:
	\begin{equation*}\label{wsubdifrad}
	f^r(\overline{x};x) \geq  \langle v,x \rangle -c\left\| x\right\| \mbox{ for all } x \in \mathbb{R}^n.
	\end{equation*}
	
	In this proof, we aim not only to construct a weak subgradient $(v,c)$ of the given function $f,$ but also to construct a maximal weak subgradient for a given point $h,$ in the sense that the function $g(x)=  \langle v,x \rangle -c\left\| x \right\|$ is
	everywhere less than or equal to $ g(h)=f^r(\overline{x};h)- \varepsilon$ and that $g$ achieves its maximum value on the unit sphere $S_1=\{x \in \mathbb{R}^n : \|x\|=1 \}$ at the point $x=h.$ 
	
	We will seek a pair $(v,c)$ such that, the directional derivative 
	$$g^{\prime}(h;y)= \langle v,y \rangle -c \langle h,y \rangle /
	\left\| h\right\| \quad \mbox{for all } y \in \mathbb{R}^n
	$$ 
	of function $g$ at $h$ in direction $y,$ equals zero on the subspace
	$$
	\mathbb{H}=
	\{y \in \mathbb{R}^n : \langle h,y \rangle=0\}.
	$$ 
	Then, the equality $g^{\prime}(h;y)= 0$ on the subspace $\mathbb{H},$ implies:
	\begin{equation*}\label{eq4.6}
	\langle v,y \rangle=0 \quad \mbox{ for all } y \in \mathbb{H}.
	\end{equation*}
	Thus we obtain that the vector $v$ must be orthogonal to the subspace $\mathbb{H}$. Since
	$\mathbb{H}$ is an $(n-1)$-dimensional subspace of $\mathbb{R}^n$, there exists a set of orthonormal
	basis vectors $\{e_1, \ldots ,e_{n-1}\}$ in $\mathbb{H}$. Then, by orthogonality of $v$ to the
	subspace $\mathbb{H},$ we have
	\begin{equation}\label{vejeqzero}
	\langle v,e_j \rangle =0 \quad \mbox{ for all } j=1,\ldots,n-1.
	\end{equation}
	Now note that the condition $g(h)=f^r(\overline{x};h)- \varepsilon$ leads to the relation
	$ \langle v,h\rangle -c\|h\|=f^r(\overline{x};h)- \varepsilon$. By using the equality $\|h\|=1$ and combining
	this equality with the $n-1$ relations given in \eqref{vejeqzero}, we obtain $n$ equations for
	$n+1$ unknown parameters $(v,c)\in \mathbb{R}^n\times \mathbb{R}_+$ in the following form:
	\begin{eqnarray}
	\langle v,h\rangle &=& c\|h\| + f^r(\overline{x};h) - \varepsilon,\label{eqn1}\\
	\langle v,e_j \rangle &=& 0 \quad \mbox{ for all }
	j=1,\ldots,n-1.\label{eqn2}
	\end{eqnarray}
	Since the vector $h$ is chosen to be perpendicular to the subspace $\mathbb{H},$ and the basis
	vectors $e_j,\quad j=1,\ldots,n-1 $ are orthonormal, we obtain that the vectors
	$h,e_1,\ldots,e_{n-1}$ are linearly independent, and therefore the system of linear equations
	given by relations $(\ref{eqn1})$--$(\ref{eqn2})$ has a unique solution $v$ for each $c$.
	
	We now find a solution to the system of equations (\ref{eqn1})--(\ref{eqn2}) explicitly.
	Recall that the vector $h$ is orthogonal to the subspace  $\mathbb{H}$. Therefore we can seek
	a solution to the set of equations (\ref{eqn2}) in the form $v=\lambda h$, where
	$\lambda$ is an unknown coefficient. By substituting this expression for $v$ in
	(\ref{eqn1}), we obtain $\lambda=c+f^r(\overline{x};h)- \varepsilon$.
	
	Thus for any given $c \geq 0,$  we have obtained a pair $(v_c,c)\in \mathbb{R}^n\times \mathbb{R}_+$ with 
	$$v_c=(c+f^r(\overline{x};h)- \varepsilon)h$$ 
	such that
	$$ g(x) \leq g(h) = f^r(\overline{x};h)- \varepsilon.$$
	
	Now we show that the number $c$ in the definition of $g,$ can be chosen large
	enough such that
	\begin{equation}\label{gx7}
	g_c(x) = \langle v_{c},x \rangle - c\left\| x\right\|   \leq f^r(\overline{x};x) \mbox { for all } x \in \mathbb{R}^n.
	\end{equation}
	
	For this aim, since $g$ and $f^r(\overline{x};\cdot)$ are both positively homogeneous functions, it is sufficient
	to show \eqref{gx7} only for points $x$ in the unit sphere $S_1.$
	
	Suppose to the contrary that there exist sequences $\{c_k\}$ with $c_k \to +\infty$ and
	$\{x_k\} \subset S_1$ such that
	\[
	g_{c_k}(x_{k}) = \langle v_{c_k},x_{k} \rangle - {c_k}\left\| x_{k}\right\| > f^r(\overline{x};x_{k}) \mbox { for all } k = 1, 2, \ldots
	\]
	or, since $\left\| x_{k}\right\| = 1,$
	\begin{equation}
	c_k (\langle h,x_k \rangle -1) + f^r(\overline{x};h)\langle h,x_k \rangle - f^r(\overline{x};x_{k}) - \varepsilon \langle
	h,x_k \rangle > 0  \label{gx8}
	\end{equation}
	for all  $k = 1, 2, \ldots$ .
	Without loss of generality we can assume that $x_k$ is a convergent sequence. Consider two
	cases.
	
	(Case 1) Let $x_k \to \widetilde{x} \neq h.$ In this case, since both $h$ and $\widetilde{x}$
	are in a unit circle, we have $\langle h,\widetilde{x} \rangle -1 < 0.$ Then, due to the
	boundedness from below of $f^r(\overline{x};\cdot)$ on the unit sphere (by the hypothesis, $f^r(\overline{x};\cdot)$ is given by \eqref{radial20}), the relation (\ref{gx8})
	leads to a contradiction for $k \to \infty.$
	
	(Case 2) Let $x_k \to  h.$ Now, since $\langle h,h \rangle =1,$  by letting to the limit as
	$k \to \infty$, we obtain $- \varepsilon >0$, which is a contradiction.
	
	Thus (\ref{gx7}) is proved, and it is shown that given any $\varepsilon > 0,$ there exists a
	number $c_\varepsilon > 0$ such that the function $g_{c_\varepsilon}$ corresponding to the
	pair $ (v_\varepsilon, c_\varepsilon) = ((c_\varepsilon + f^r(\overline{x};h) - \varepsilon)h, c_\varepsilon)$,
	defined as
	$$
	g_{c_\varepsilon} (x) = (c_\varepsilon + f^r(\overline{x};h) - \varepsilon) \langle h,x \rangle -c_\varepsilon \left\|
	x \right\|
	$$
	satisfies the following conditions
	$$
	g_{c_\varepsilon} (x) \leq f^r(\overline{x};x) \mbox{  for all  } x \in \mathbb{R}^n,
	$$
	and
	$$
	g_{c_\varepsilon} (x) \leq g_{c_\varepsilon} (h) = f^r(\overline{x};h) - \varepsilon.
	$$
	The first relation, in particular, means that $ (v_\varepsilon, c_\varepsilon) \in \partial^wf^r(\overline{x};0).$ Hence by Lemma \ref{weaksubdifrad}, we obtain that, $ (v_\varepsilon, c_\varepsilon) \in \partial^w f(\overline{x}),$ which completes the proof.
	$\Box$
\end{proof}



Now we give the $\ell_1-$ norm version of Theorem \ref{computingvall2}. Since the proof of this version, is similar to that of Theorem \ref{computingvall2}, we present it without the proof.


\begin{theorem} \label{computingvall1}
	Let $f:\mathbb{R}^n\rightarrow \mathbb{R}$ be a proper  function, $\overline{x} \in (\mathbb{R}^n,\|\cdot\|_1)$,
	and $\overline{y}=f(\overline{x})$ be finite. Assume that, $f$ is radially epidifferentiable at $\overline{x}.$ Then, for every $\varepsilon >0$ and $h =(h_1,\ldots,h_n) \in \mathbb{R}^n \setminus \{0_{\mathbb{R}^n}\},$ there exists a weak subgradient  $(v,c)  \in \partial^wf(\overline{x})$ such that
	\begin{equation}\label{weaksubdircomp1}
	v=(c+f^r(\overline{x};h)- \varepsilon)Sgn(h)
	\end{equation}
	or 
	\begin{equation}\label{radialepidercomp1}
	f^r(\overline{x};h) = \frac{1}{n} \langle v,Sgn(h) \rangle - c + \varepsilon,
	\end{equation}
	where $Sgn(h)$ is the $n$ dimensional vector defined as $Sgn(h) = (\sgn(h_1),\sgn(h_2),\ldots,\sgn(h_n))$
	and $\sgn(h_i) = 1 $ if $h_i >0$ and $\sgn(h_i) = -1 $ if $h_i <0, i=1,\ldots,n.$ 
\end{theorem}


\begin{example}\label{ex3weak}
	Consider the function $f_3$ from Example \ref{example3} (see Figure \ref{figure:f_3_RE}) and try to illustrate Theorems \ref{computingvall2} and \ref{computingvall1}. Recall that
	
	\begin{equation*}
	f_3(x) =
	\left\{
	\begin{array}{ll}
	4|x+1|    & \mbox{ if } x \leq 0, \\
	|x-1| + 3 & \mbox{ if } x > 0.
	\end{array}
	\right.
	\end{equation*}
	
	For this function, we will compute the weak sudifferentials at different points.\\
	
	First consider the point $\overline{x} = 1.$
	By Lemma \ref{weaksubdifrad} we have $\partial^{w}f^{r}_3(1;0) = \partial^{w} f_3(1).$ By definition of the weak subdifferential we obtain:	
	\begin{eqnarray} \label{wsubdiff3at1}
	& & \partial^{w} f_3(1) =	\partial^{w}f_3^r(1;0) \nonumber \\
	&=& \{(v,c) \in \mathbb{R} \times \mathbb{R}_+ : f_3^r(1;h) - f_3^r(1;0) \geq vh - c|h| \mbox{ for all } h \in \mathbb{R}  \} \nonumber\\
	&=& \{(v,c) \in \mathbb{R} \times \mathbb{R}_+ : -c-\frac{3}{2} \leq v \leq c+1 \}.
	\end{eqnarray}

	Now try to compute weak subgradients by applying Theorems \ref{computingvall2} and/or \ref{computingvall1}. 
	Let $h=-1, \varepsilon=1/2.$ Then by \eqref{radepiderex3} we have $f^r_3(1;-1)=-3/2,$ and applying formula $v=(c+f^r_3(\overline{x};h)- \varepsilon)h$ we obtain $v= -c+2.$ By checking with \eqref{wsubdiff3at1}, we see that $(v,c) = (-c+2, c) \in \partial^{w} f_3(1)$ for every $c \geq 1/2.$ 
	
	$\Box$
\end{example}


The following section discusses and studies optimality conditions via the radial epiderivatives.



\section{Optimality conditions via generalized derivatives}\label{optcond}

We begin this section with the following result which establishes a necessary and sufficient condition for a descent direction via the radial epiderivative, for nonconvex nonsmooth functions. We will say that $h \in \mathbb{X}$ is a descent direction for function $f:\mathbb{X} \rightarrow \mathbb{R} \cup \{+\infty\}$ at $\overline{x} \in \mathbb{X},$ if there exists a pozitive number $\bar{t}$ such that $f(\overline{x}+\bar{t}h)<f(\overline{x}).$

 
\begin{theorem} \label{radepiderdescent}
	Let $(\mathbb{X}, \|.\|_\mathbb{X})$ be real normed space and let $f:\mathbb{X} \rightarrow \mathbb{R} \cup \{+\infty\}$  be a proper function. Assume that $f$ is radially epidifferentiable at $\overline{x} \in \mathbb{X}.$ Then the vector $h \in \mathbb{X}$ is a descent direction for $f$ at $\overline{x}$ if and only if   $f^r(\overline{x}; h) < 0.$  
\end{theorem}	
\begin{proof} \textbf{Proof of If.}
	let $f$ be radially epidifferentiable at $\overline{x}.$ Assume that  $f^r(\overline{x}; h) < 0$ for some $h \in  \mathbb{X}.$ Then by Proposition \ref{rews} we have:
	\begin{equation*} 
	f^{r}(\overline{x};h) = \inf_{t > 0} \lim_{ u \rightarrow h} \frac{f(\overline{x}+ tu)- f(\overline{x})}{t} <0.
	\end{equation*}
	Then, there exists a pozitive number $\varepsilon >0$ such that $f^{r}(\overline{x};h) < -\varepsilon.$ Thus there exists a pozitive number $t_{\varepsilon}$ such that
	\begin{equation*} 
	\lim_{ u \rightarrow h} \frac{f(\overline{x}+ t_{\varepsilon}u)- f(\overline{x})}{t_{\varepsilon}} < -\varepsilon.
	\end{equation*}
	Since $f$ is radially epidifferentiable at $\overline{x},$ by Theorem \ref{radepiderlsc} it is lower semicontinuous there and hence, the latter relation implies:
	\begin{equation*} 
	\frac{f(\overline{x}+ t_{\varepsilon}h)- f(\overline{x})}{t_{\varepsilon}} < -\varepsilon,
	\end{equation*}
	which means that  $h$ is a descent direction for $f$ at $\overline{x}.$
	
	The proof of \textbf{``only if"} is similar to that of \textbf{``if"} part.
	
	$ \Box $ 
\end{proof}


Now we formulate the following optimality condition which can easily be obtained from Theorem \ref{radepiderdescent}. Note that the similar optimality condition was earlier established by Kasimbeyli in \cite[Theorem 3.6]{Kasimbeyli2009}.

\begin{corollary} \label{radepiderglobalmin}
	Let $(\mathbb{X}, \|.\|_\mathbb{X})$ be real normed space and let $f:\mathbb{X} \rightarrow \mathbb{R} \cup \{+\infty\}$  be a proper function. Assume that $f$ is radially epidifferentiable at $\overline{x} \in \mathbb{X}.$ Then $f$ attains global minimum at  $\overline{x}$ if and only if   $f^r(\overline{x}; h)$ attains its minimum at $h=0.$  
\end{corollary}	
\begin{proof} The proof easily follows from Theorem \ref{radepiderdescent}.
	
	$ \Box $ 
\end{proof}


\begin{example}\label{ex13}
	Consider the function $f_3$ from Example \ref{example3}. Let
	
	\begin{equation*}
	f_3(x) =
	\left\{
	\begin{array}{ll}
	4|x+1|    & \mbox{ if } x \leq 0, \\
	|x-1| + 3 & \mbox{ if } x > 0.
	\end{array}
	\right.
	\end{equation*}
	
	Obviously $h=0$ is a global minimum of the radial epiderivative 
	$$
	f^r_3(-1;h) = 	\left\{
	\begin{array}{ll}
	-4h    & \mbox{ if } h \leq 0, \\
	h      & \mbox{ if } h > 0
	\end{array}
	\right.
	$$ 
	(see \eqref{radepiderex3}), which illustrates the assertion of Corollary \ref{radepiderglobalmin} and demonstrates that $\overline{x} =-1$ is a global minimum of $f_3.$ On the other hand, since  $\overline{x} =-1$ is a global minimum of $f_3,$ at this point the other (global) optimality condition $(0,0) \in \partial^{w}f_3(-1)$ must be satisfied (see \eqref{optcondzeroweaksubgrad}). It can easily be checked that $(0,0) \in \partial^{w}f_3(-1) = \{(v,c) \in \mathbb{R} \times \mathbb{R}_+ : -c-4 \leq v \leq c+1 \}.$  
	
	As another illustration of Corollary \ref{radepiderglobalmin}, consider the point $\overline{x} =1,$ which is a local (but not global) minimum of $f_3.$ At this point (see \eqref{radepiderex3}) we have:
	$$ 
	f^r_3(1;h)= \left\{
	\begin{array}{ll}
	\frac{3h}{2}  & \mbox{ if } h \leq 0, \\
	h    & \mbox{ if } h > 0. \\
	\end{array}
	\right.	
	$$
	Clearly, $h=0$ is not a minimum point of $f^r_3(1;h),$ and as a result, it can easily be seen that $(0,0) \notin \partial^{w}f_3(1).$
	
	On the other hand, since $f^r_3(1;h) <0$ for every $h<0,$ we obtain that (for example) $h=-1 $ is a descent direction for $f_3$ at $\overline{x} = 1.$ In this case, the better point can be computed in the form  $x = 1 + t (-1)$ and the optimal value for $t=t_{opt}>0$ can be found by solving the scalar problem : $\min \{ f_3(1-t) : t>0 \}.$ An easy computing shows that for $t=2 $ the next iteration gives the global minimum $x = -1.$
	
	$\Box$
\end{example}


\begin{remark}
	Note that Theorem \ref{radepiderdescent} and Corollary \ref{radepiderglobalmin} allow one to determine a descent direction even in the case when $\overline{x}$ is a local minimum of the (nonconvex) function under consideration. This means that by using this theorem, one can actually devise a minimization method for finding a global minimum. 
\end{remark}


\begin{remark}
	In \cite{DincyalcinK2020opt} the authors developed a method for approximate computing the weak subgradients via the directional derivatives, which was used there to formulate a solution algorithm for solving some classes of nonconvex optimization problems. With the help of Theorems \ref{computingvall2} or \ref{computingvall1}, the approximate computing method for weak subgradients can be used to estimate radial epiderivatives and by this way,  using Theorem \ref{radepiderdescent} one can compute (global) descent direction for a (nonconvex) function under consideration.
	On the other hand, if we are given the value of the radial epiderivative, we can estimate the weak subgradients by using Theorems \ref{computingvall2} or \ref{computingvall1}, and then use them in the weak subgradient based solution method given in \cite{DincyalcinK2020opt}.
\end{remark}


Finally we illustrate the behaviour of the generalized derivatives considered in this paper, and optimality conditions given in Theorem \ref{clarke-stationary-point} and in relations  \eqref{optcondzeroweaksubgrad} and \eqref{optcondzeroradepider} on two simple functions.

\begin{example} \label{exsubdif}
	Let
	\begin{equation}
	f_8(x) =
	\left\{
	\begin{array}{ll}
	x^2    & \mbox{ if } x \leq 0, \\
	-x + 1 & \mbox{ if } x > 0.
	\end{array}
	\right.
	\end{equation}
	
	\begin{figure}[h]
		\scalebox{0.4}{\includegraphics{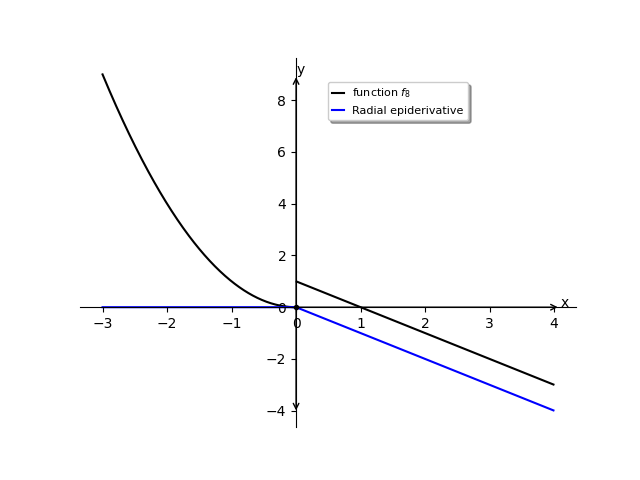}}
		\scalebox{0.4}{\includegraphics{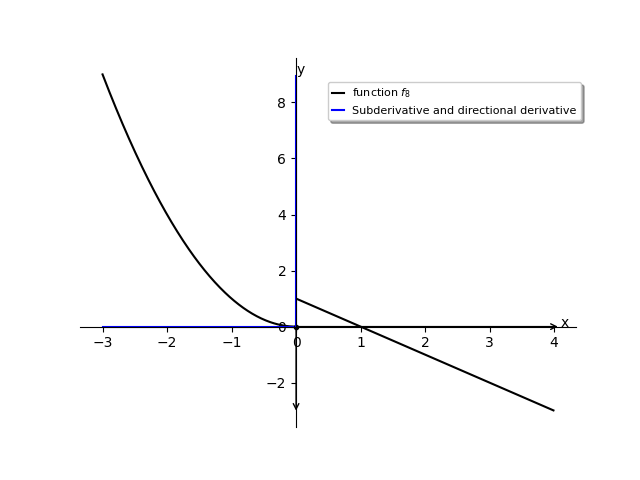}}
		\caption{The radial epiderivative (left), and subderivative and Clarke's derivative (right) of $f_8$ at the point $\overline{x}=0$}
		\label{figure:f_1}
	\end{figure}
	Then, 
	$$
	f_8^r(0;h)= \left\{
	\begin{array}{ll}
	0    & \mbox{ if } h \leq 0, \\
	-h   & \mbox{ if } h > 0.
	\end{array}
	\right.
	$$
	Clearly, $f_8^r(0;h)\rightarrow -\infty$ as $h \rightarrow +\infty$ which indicates that function $f_8$ is unbounded from below and hence has no a global minimum value, see Figure \ref{figure:f_1}.
	
	On the other hand,
	$$
	\text{d}f_8(0;h) = \left\{
	\begin{array}{ll}
	0    & \mbox{ if } x \leq 0, \\
	+\infty   & \mbox{ if } x > 0.
	\end{array}
	\right.
	$$
	Consequently, $f_8^{\circ}(0;h)= \text{d}f_8(0;h).$ Despite the fact that $x=0$ is not a minimum point of $f_8,$ the generalized subdifferential set contains the zero element, indicating that the point $x=0$ is a Clarke stationary point. 
	
	Now consider the function
	\begin{equation}
	f_9(x) =
	\left\{
	\begin{array}{ll}
	x^2    & \mbox{ if } x \leq 0, \\
	x + 1 & \mbox{ if } x > 0.
	\end{array}
	\right.
	\end{equation}
	
	\begin{figure}[h]
		\scalebox{0.4}{\includegraphics{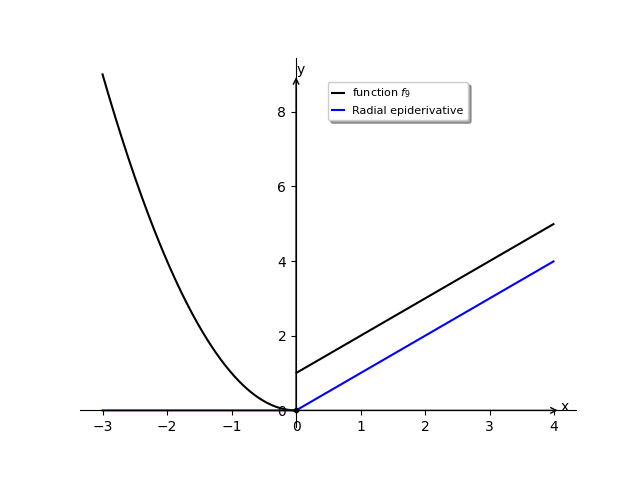}}
		\scalebox{0.4}{\includegraphics{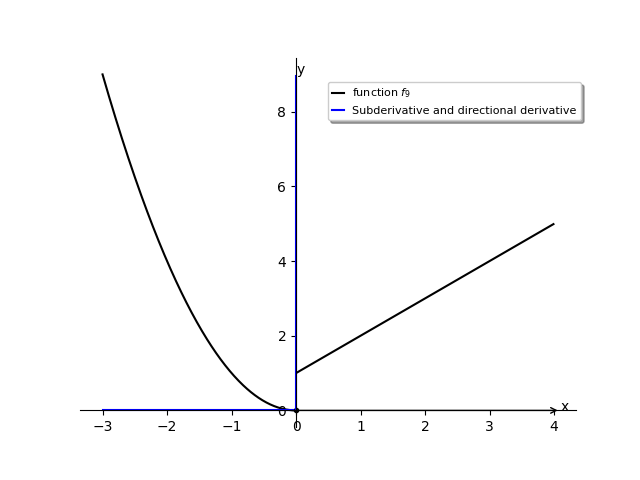}}
		\caption{The radial epiderivative (left), and subderivative and Clarke's derivative (right) of $f_9$ at the point $\overline{x}=0$}
		\label{figure:f_2}
	\end{figure}

	Then, 
	$$
	f_9^r(0;h)= \left\{
	\begin{array}{ll}
	0    & \mbox{ if } x \leq 0, \\
	h   & \mbox{ if } x > 0.
	\end{array}
	\right.
	$$
	Clearly, $f_9^r(0;h)$ attains its global minimum value at $h=0$ indicating that the function $f_9$ attains its global minimum at $\overline{x}=0,$ see Figure \ref{figure:f_2}. Note also that $(0,0) \in \partial^w f_9(0)$ which again justifies this assertion.
	
	On the other hand, despite the differences between the functions $f_8$ and $f_9$ we obtain the same expressions for the Clarke directional derivative and the subderivative: $\text{d}f_8(0;h) = \text{d}f_9(0;h)=f_8^{\circ}(0;h)=f_9^{\circ}(0;h)=+\infty$ for $h>0.$
	$\Box$
\end{example}	 



\section{Conclusion}\label{conclusion} This paper studies new properties of the directional derivatives, the Rockafellar's subderivatives, the Clarke derivatives and the radial epiderivatives. These properties are analyzed and illustrated on examples. The paper presents new regularity conditions for establishing equality between the generalized derivatives and compares them with the existing in the literature conditions. We present a new formula for defining the radial epiderivatives and establish new optimality conditions via the radial epiderivatives. These conditions are compared with the optimality conditions given in the literature via the generalized derivatives. All the regularity and optimality conditions are demonstrated and illustrated on illustrative examples. By using the global optimality condition studied in the paper, we explain  an idea for developing a solution method in nonconvex optimization, which can be a subject for the new research directions in the future. The paper also presents explicit formulas for approximate computing the weak subgradients in terms of the radial epiderivatives and vice versa, which again can be used to develop a new radial epiderivatives based global solution method in nonconvex programming.


\section*{Statements and Declarations}

\textbf{Competing Interests and Funding.} This study is supported by The Scientific and Technological Research Council of Turkey (TUBITAK) under the Grant No 217M487. The authors declare that they have no any conflict of interests.






%
%


%
%

\end{document}